\documentclass[12pt]{article}
\usepackage{amssymb}
\newtheorem{theorem}{Theorem}

\newtheorem{proposition}[theorem]{Proposition}

\newtheorem{lemma}[theorem]{Lemma}
\newtheorem{corollary}[theorem]{Corollary}
\newtheorem{remark}[theorem]{Remark}

\newtheorem{example}[theorem]{Example}

\newcommand{\R}{\mathbb{R}}
\newcommand{\Q}{\mathbb{Q}}

\newcommand{\Sf}{\mathbb{S}}

\newcommand{\Hy}{\mathbb{H}}

\newcommand{\spa}{\mbox{span}}

\newcommand{\Ima}{\mbox{Im }}

\newcommand{\rank}{\mbox{rank }}

\newcommand{\po}{{\hspace*{-1ex}}{\bf .  }}

\def\<{\langle}
\def\n{\nabla}
\def\>{\rangle}
\def\a{\alpha}

\def\e{\epsilon}

\def\bea{\begin{eqnarray*} }
\def\eea{\end{eqnarray*} }
\def\be{\begin{equation} }
\def\ee{\end{equation} }

\def\proof{\noindent{\it Proof: }}
\def\qed{\ifhmode\unskip\nobreak\fi\ifmmode\ifinner
\else\hskip5 pt \fi\fi\hbox{\hskip5 pt \vrule width4 pt
height6 pt  depth1.5 pt \hskip 1pt }}
\begin{document}
\title{\bigskip
\bigskip
Submanifolds of  products of  space forms.
}
\author{ B. Mendon\c ca and R. Tojeiro}
\date{}
\maketitle

\begin{abstract}
We give a complete classification of submanifolds with parallel second fundamental form of a 
product of two space forms. We also reduce the classification of  umbilical submanifolds with dimension $m\geq 3$ of a  product $\Q_{k_1}^{n_1}\times \Q_{k_2}^{n_2}$ of two space forms whose  curvatures satisfy $k_1+k_2\neq 0$ to the classification of $m$-dimensional umbilical submanifolds of codimension  two of $\Sf^n\times \R$ and $\Hy^n\times \R$. The case of  $\Sf^n\times \R$ was carried out in \cite{mt}. As a main tool  we derive reduction of codimension theorems of independent interest for  submanifolds of products of two space forms.

\end{abstract}

\vspace{2ex}


\section[Introduction]{Introduction}

Let $f\colon\, M\to N$ be an isometric immersion between Riemannian manifolds and let $\alpha\colon\,TM\times TM\to N_fM$
be its second fundamental form with values in the normal bundle. Then $f$ is said to have \emph{parallel second fundamental form} if 
$\nabla \alpha=0$, where $\nabla \alpha$ is the  Van der Waerden-Bortolloti covariant derivative of $\alpha$. One also says for short that \emph{ $f$ is parallel}. Roughly speaking, this means that $f$ has the same second fundamental form at any point of $M$, if tangent and normal spaces at any two distinct points are identified by means of  parallel transport in the tangent and normal connections, respectively,  along any curve  joining them.

Parallel submanifolds of Euclidean space have been classified by Ferus \cite{fe}. He showed that all of them are products of an Euclidean factor and standard minimal embeddings into hyperspheres of symmetric $R$-spaces, which are orbits of special types of  $s$-representations. The case of the sphere is an easy consequence of the Euclidean one, whereas the classification of parallel submanifolds of hyperbolic space was carried out independently by  Backes--Reckziegel \cite{br} and Takeuchi  \cite{ta}.

Apart from space forms, however, parallel submanifolds of a Riemannian manifold have been classified only in a few other cases, e.g. for 
simply connected rank one symmetric spaces (see, e.g., the discussion in Chapter $9$ of \cite{bco}).  

One of our main results is a complete classification of  parallel submanifolds  of a 
product of two space forms. We state separately the cases in which one of the factors is flat or not. 
First  observe that, given  $k_1, k_2\in \R$ with $k_1k_2>0$, the map 
\be \label{eq:g}g\colon\,\Q_k^n\to \Q_{k_1}^{n}\times  \Q_{k_2}^{n}, \,\,\,\,g(x)=(ax, bx),\ee where $k=k_1k_2/(k_1+k_2)$,  $a^2=k_2/(k_1+k_2)$ and  $b^2=k_1/(k_1+k_2)$, is a totally geodesic isometric embedding  (see Example \ref{ex:tgeod} below).


We say that $f\colon\,M^m\to \Q_{k_1}^{n_1}\times \Q_{k_2}^{n_2}$ is an isometric immersion into a slice of $\Q_{k_1}^{n_1}\times \Q_{k_2}^{n_2}$ if there exist an isometric immersion $\tilde f\colon\,  M^m\to \Q_{k_1}^{n_1}$ (resp., $\tilde f\colon\,  M^m\to \Q_{k_2}^{n_2}$) and a totally geodesic inclusion $j\colon\, \Q_{k_1}^{n_1}\to \Q_{k_1}^{n_1}\times \{x_2\}$ (resp., $j\colon\, \Q_{k_2}^{n_2}\to \{x_1\}\times \Q_{k_2}^{n_2}$) such that $f=j\circ \tilde f$.

\begin{theorem}\po\label{thm:parallel} Let $f\colon\,M^m\to \Q_{k_1}^{n_1}\times \Q_{k_2}^{n_2}$, $k_1k_2 \neq 0$,  be a parallel isometric immersion. Then one of the following possibilities holds:
\begin{itemize}
\item[$(i)$]  $f$ is a parallel isometric immersion into a slice of $\Q_{k_1}^{n_1}\times \Q_{k_2}^{n_2}$;
\item[$(ii)$] $M^m$ is locally a Riemannian product $M^m=M_1^{m_1}\times M_2^{m_2}$ and $f=f_1\times f_2$, where $f_i\colon\,M^{m_i}_i\to \Q_{k_i}^{n_i}$, $1\leq i\leq 2$, is a parallel isometric immersion. 
\item[$(iii)$] $k_1k_2>0$ and there exists a parallel isometric immersion $\bar f\colon M^m\to \Q_{k}^{m+\ell}$, with $k=k_1k_2/(k_1+k_2)$,  such that
$f=j\circ g\circ \bar f$, where $j\colon\, \Q_{k_1}^{m+\ell}\times \Q_{k_2}^{m+\ell}\to \Q_{k_1}^{n_1}\times \Q_{k_2}^{n_2}$ is a totally geodesic inclusion and $g\colon\, \Q_{k}^{m+\ell}\to \Q_{k_1}^{m+\ell}\times \Q_{k_2}^{m+\ell}$ is as in (\ref{eq:g}).
\end{itemize}
Moreover, the second possibility holds globally if $M^n$ is complete and simply connected.
\end{theorem}

 Notice that  the last possibility can occur  only if both $n_1$ and $n_2$ are greater than or equal to  $m$. Moreover,  if  this is the case and either $n_1$ or $n_2$ is equal to $m$ then $\ell=0$ and  $\bar f$ is just a local isometry, in which case $f$ is totally geodesic. 
 
 We point out that, after a preliminary version of this article was completed, we learned that the case $k_1=k_2\neq 0$ of Theorem \ref{thm:parallel} was independently obtained with a different approach by Jentsch \cite{je}  as a consequence of his classification of parallel submanifolds of the Grassmannian $G_2^+(\R^{n+2})$ of oriented $2$-planes of $\R^{n+2}$ and of its noncompact dual.\vspace{1ex}
 
 In case one of the factors is flat, the classification of parallel submanifolds of a product of space forms reads as follows. We recall that a parallel  unit speed curve $\gamma\colon\, \R\to M$ on a Riemannian manifold is also called an \emph{extrinsic circle}. Thus, $\gamma$ is an extrinsic circle  if its curvature vector $\nabla_{\gamma'}\gamma'$ is parallel with respect to its normal connection.  
By an extrinsic circle $\gamma\colon\,\R\to \Q_{k_1}^{2}\times \R$ being \emph{full} we mean that $\gamma(\R)$ does not lie in a totally geodesic surface of $\Q_{k_1}^{2}\times \R$.

 \begin{theorem}\po\label{thm:parallel2} Let $f\colon\,M^m\to \Q_{k_1}^{n_1}\times \R^{n_2}$, $k_1\neq  0$,  be a parallel isometric immersion. Then one of the following possibilities holds:
\begin{itemize}
\item[$(i)$]  $f$ is a parallel isometric immersion into a slice of $\Q_{k_1}^{n_1}\times \R^{n_2}$;
\item[$(ii)$] $M^m$ is locally a Riemannian product $M^m=M_1^{m_1}\times M_2^{m_2}$ and $f=f_1\times f_2$, where $f_1\colon\,M^{m_1}_1\to \Q_{k_1}^{n_1}$ and $f_2\colon\,M^{m_2}_1\to \R^{n_2}$ are parallel isometric immersions.
\item[$(iii)$]  $k_1>0$ (resp., $k_1<0$) and $f\circ \tilde \Pi=j\circ  \Pi\circ \tilde f$ (resp., $f=j\circ  \Pi\circ \tilde f$), where $\tilde \Pi\colon\, \tilde M^m\to M^m$ is the universal covering of $M^m$, $\tilde f\colon\, \tilde M^m\to \R^{n_2+1}$ (resp., $\tilde f\colon\,  M^m\to \R^{n_2+1}$)  is  a parallel isometric immersion,  $j\colon\, \Q_{k_1}^{1}\times  \R^{n_2}\to \Q_{k_1}^{n_1}\times \R^{n_2}$ is  totally geodesic  and $ \Pi\colon\, \R^{n_2+1}\to\Q_{k_1}^{1}\times  \R^{n_2}$ is  a locally isometric covering map (resp., isometry).
\item[$(iv)$] $M^m$ is locally a Riemannian product $M^m=\R\times N^{m-1}$ and $f=j\circ(\gamma\times \tilde f)$, where $j\colon\, \Q_{k_1}^{2}\times  \R^{n_2}\to \Q_{k_1}^{n_1}\times \R^{n_2}$ is a totally geodesic inclusion, $\gamma\colon\,\R\to \Q_{k_1}^{2}\times \R$ is a full extrinsic circle and $\tilde f\colon\,N^{m-1}\to \R^{n_2-1}$ is a parallel isometric immersion. 

\end{itemize}
Moreover, the second and fourth possibilities hold globally if $M^n$ is complete and simply connected.
\end{theorem}

Notice that case  $(iv)$ (respectively, $(iii)$) can occur only if $n_2\geq m$ (respectively, $n_2\geq m-1$), and $f$ must be totally geodesic if equality holds.  
Therefore, a nontotally geodesic parallel isometric immersion $f\colon\,M^m\to \Q_{k_1}^{n_1}\times \R^{n_2}$, $k_1\neq  0$, $n_2\leq m-1$, must be either as in  $(ii)$ or a  parallel isometric immersion into a slice $\Q_{k_1}^{n_1}\times \{x_2\}\subset \Q_{k_1}^{n_1}\times \R^{n_2}$. Moreover, in the last case one must have $n_1\geq m+1$. This extends the results in \cite{ckv} and \cite{vv}  for the case of hypersurfaces of $\Q_{k}^{n}\times \R$. 

As a consequence of Theorems \ref{thm:parallel} and \ref{thm:parallel2}, we obtain the  classification of totally geodesic submanifolds of $\Q_{k_1}^{n_1}\times \Q_{k_2}^{n_2}$.

\begin{corollary}\po\label{cor:totgeodesic} Let $f\colon\,M^m\to \Q_{k_1}^{n_1}\times \Q_{k_2}^{n_2}$, $(k_1,k_2) \neq (0,0)$,  be a totally geodesic isometric immersion. Then one of the following possibilities holds:
\begin{itemize}
\item[$(i)$]  $f$ is a totally geodesic isometric immersion into a slice of $\Q_{k_1}^{n_1}\times \Q_{k_2}^{n_2}$;
\item[$(ii)$] There exist a local isometry $\phi\colon\,  M^m\to \Q_{k_1}^{m_1}\times \Q_{k_2}^{m_2}$ and a totally geodesic inclusion $j=j_1\times j_2\colon\, \Q_{k_1}^{m_1}\times \Q_{k_2}^{m_2}\to \Q_{k_1}^{n_1}\times \Q_{k_2}^{n_2}$ such that $f=j\circ \phi$.
\item[$(iii)$] $k_1k_2>0$ and there exist a local isometry $\phi\colon\,  M^m\to \Q_{k}^{m}$, $k=k_1k_2/(k_1+k_2)$, and a totally geodesic inclusion $j=j_1\times j_2\colon\, \Q_{k_1}^{m}\times \Q_{k_2}^{m}\to \Q_{k_1}^{n_1}\times \Q_{k_2}^{n_2}$ such that $f=j\circ g\circ \phi$, where $g\colon\, \Q_{k}^{m}\to \Q_{k_1}^{m}\times \Q_{k_2}^{m}$ is as in (\ref{eq:g}).
\item[$(iv)$] $k_1k_2=0$, say, $k_2=0$, and there exist a local isometry $\phi\colon\,  M^m\to \R^m=\R\times \R^{m-1}$, a unit-speed geodesic $\gamma\colon\, \R\to \Q^1_{k_1}\times \R$ and a totally geodesic inclusion $j\colon\, \Q_{k_1}^{1}\times \R^{m}\to \Q_{k_1}^{n_1}\times \R^{n_2}$ such that $f=j\circ (\gamma\times id)\circ \phi$, where $id\colon\, \R^{m-1}\to \R^{m-1}$ is the identity.
\end{itemize}
\end{corollary}

Throughout the paper we use the framework introduced in \cite{ltv} for studying submanifolds of products of space forms.
As a main tool for the proofs of Theorems \ref{thm:parallel} and~\ref{thm:parallel2} we derive reduction of codimension theorems
of independent interest for arbitrary submanifolds of products of space forms, some of which extend well-known results for submanifolds of space forms. 

In the last part of the paper we apply them to
study umbilical submanifolds of a product of two space forms. Recall that an isometric immersion $f\colon\, M^m\to \tilde{M}^n$ between Riemannian manifolds  is  \emph{umbilical}
 if there exists a normal vector field $\zeta$ along $f$ such that its second fundamental form  satisfies $\alpha(X,Y)=\<X,Y\>\zeta$ for all $X,Y\in TM$. One main motivation for this study is a theorem by Nikolayevsky (see Theorem $1$ of \cite{al}), which states that any umbilical submanifold of a symmetric space $N$ is an umbilical submanifold of a product of space forms totally geodesically embedded in $N$.
 
 We prove the following result, which reduces the classification of $m$-dimensional umbilical submanifolds of $\Q_{k_1}^{n_1}\times \Q_{k_2}^{n_2}$, $m\geq 3$ and $k_1+k_2\neq 0$, to the cases in which either $n_1\in \{m, m+1\}$ and $n_2=1$ or $n_2\in \{m, m+1\}$ and $n_1=1$, or equivalently, by passing to the universal coverings, to the classification of $m$-dimensional umbilical submanifolds of codimension two of $\Sf^n\times \R$ and $\Hy^n\times \R$. Here $\Sf^n$ and $\Hy^n$ stand for the sphere and hyperbolic space, respectively. The case of  $\Sf^n\times \R$ was carried out in \cite{mt}, extending previous results in \cite{st} and \cite{vv} for  hypersurfaces.

\begin{theorem}\po\label{thm:umb} Let $f\colon\,M^m\to \Q_{k_1}^{n_1}\times \Q_{k_2}^{n_2}$, with $m\geq 3$ and $k_1+k_2\neq 0$,  be an umbilical nontotally geodesic isometric immersion. Then one of the following possibilities holds:
\begin{itemize}
\item[$(i)$]  $f$ is an umbilical isometric immersion into a slice of $\Q_{k_1}^{n_1}\times \Q_{k_2}^{n_2}$;
\item[$(ii)$] there exist umbilical isometric immersions $f_i\colon\, M^m\to \Q_{\tilde k_i}^{n_i}$, $1\leq i\leq 2$, with $\tilde{k}_1=k_1\cos^2\theta$ and $\tilde k_2=k_2\sin^2\theta$ for some $\theta\in (0, \pi/2)$,  such that $f=(\cos \theta f_1,\sin\theta f_2)$.
\item[$(iii)$] after possibly reordering the factors, we have $k_1>0$ (resp., $k_1\leq 0$) and $f\circ \tilde\Pi=j\circ  \Pi\circ \tilde f$ (resp., $f=j\circ  \Pi\circ \tilde f$), where $\tilde \Pi\colon\, \tilde M^m\to M^m$ is the universal covering of $M^m$, $\tilde f\colon\,\tilde M^m\to \R\times \Q_{k_2}^{m+\delta}$ (resp., $\tilde f\colon\, M^m\to \R\times \Q_{k_2}^{m+\delta}$)  is  an umbilical isometric immersion with  $\delta\in \{0,1\}$, $j\colon\, \Q_{k_1}^{1}\times \Q_{k_2}^{m+\delta}\to \Q_{k_1}^{n_1}\times \Q_{k_2}^{n_2}$ is  totally geodesic  and $ \Pi\colon\,\R\times \Q^{m+\delta}_{k_2}\to  \Q_{k_1}^1\times \Q^{m+\delta}_{k_2}$ is  a locally isometric covering map (resp., isometry).
\end{itemize}
\end{theorem}

In particular, it follows from Theorem \ref{thm:umb} that there does not exist an umbilical nontotally geodesic $m$-dimensional submanifold of $\Q_{k_1}^{n_1}\times \Q_{k_2}^{n_2}$, $k_1+k_2\neq 0$, if $m$ is greater than both $n_1$ and $n_2$.

 We point out that umbilical nontotally geodesic submanifolds of $\Q_{k_1}^{n_1}\times \Q_{k_2}^{n_2}$  have been alternatively described by Nikolayevsky  as intersections of $\Q_{k_1}^{n_1}\times \Q_{k_2}^{n_2}$ with its osculating spaces at generic points in the flat underlying ambient space (see Theorem $2$ of \cite{al}).

\section[Preliminaries]{Preliminaries}

Let $\pi_i\colon\, \Q_{k_1}^{n_1}\times \Q_{k_2}^{n_2}\to \Q_{k_i}^{n_i}$ denote the canonical projection, $i=1,2$.
By abuse of notation, we denote by the same letter  its derivative, which we  regard as a  section of either $T({\Q_{k_1}^{n_1}\times \Q_{k_2}^{n_2}})^*\otimes T\Q_{k_i}^{n_i}$
or $T({\Q_{k_1}^{n_1}\times \Q_{k_2}^{n_2}})^*\otimes T({\Q_{k_1}^{n_1}\times \Q_{k_2}^{n_2}})$.
Then, the  curvature tensor
$\bar{{\cal R}}$ of $\Q_{k_1}^{n_1}\times \Q_{k_2}^{n_2}$ can be written as
$$\bar{{\cal R}}(X,Y)=k_1(X\wedge Y-X\wedge\pi_2 Y-\pi_2 X\wedge Y) +(k_1+k_2) \pi_2 X\wedge\pi_2 Y,$$
where $(X\wedge Y)Z=\<Y,Z\>X-\<X,Z\>Y$.

 Let  $f\colon\,M\to \Q_{k_1}^{n_1}\times \Q_{k_2}^{n_2}$ be  an isometric immersion of a Riemannian manifold.
  Denote by ${\cal R}$  and ${\cal R}^\perp$ the curvature tensors of the tangent and normal bundles
  $TM$ and $N_f M$, respectively,  by $\alpha=\alpha_f\in \Gamma(T^*M\otimes T^*M\otimes N_f M)$ the second fundamental form of $f$ and by $A_\eta=A^f_\eta$ its shape operator  in the normal direction $\eta$, given by  $$\<A_\eta X,Y\>=\<\alpha(X,Y),\eta\>$$ for all $X, Y\in TM$. Set
 $$L=L^{{f}}:=\pi_2\circ {f}_*\in \Gamma(T^*M\otimes T\Q_{k_2}^{n_2})\,\,\,\mbox{and}\,\,\,K=K^{{f}}:=\pi_2|_{N_fM}\in \Gamma((N_fM)^*\otimes T\Q_{k_2}^{n_2}).$$
We can write
\be\label{eq:L}
L=f_*R+S\,\,\,\,\mbox{and}\,\,\,\,K=f_*S^t+T,
\ee
where 
$$R=R^f:=L^tL,\,\,\,\,S=S^f:=K^tL\,\,\,\,\mbox{and}\,\,\,\,T=T^f:=K^tK.$$

   The tensors $R$, $S$ and $T$ were introduced in  \cite{ltv}. Note that $R$ and $T$ are symmetric. Using (\ref{eq:L}), one can check by applying $\pi_2^2=\pi_2$ to tangent and normal vectors, and then taking tangent and normal components, that they 
  satisfy the algebraic relations
\be\label{eq:pi1}
S^tS=R(I-R),\,\,\,\,\,TS=S(I-R)\,\,\,\,
\,\mbox{and}\,\,\,\,\,\,
SS^t=T(I-T).
\ee
In particular, from the first and third equations, respectively, it follows that $R$ and $T$ are in fact nonnegative operators whose eigenvalues lie in $[0,1]$.
On the other hand, taking tangent and normal components in $\nabla \pi_2=0$ and using the Gauss and Weingarten formulae  yields the differential equations
\be\label{eq:pid1}(\n_XR)Y=A_{SY}X+S^t\a(X,Y),\ee
\be\label{eq:pid2}(\n_XS)Y=T\a(X,Y)-\a(X,RY)\ee
and
\be\label{eq:pid4}(\n_X T)\eta=-SA_\eta X-\a(X,S^t\eta).\ee

 The Gauss, Codazzi and Ricci equations of $f$ are, respectively,
 \be\label{eq:gauss}
 \begin{array}{l}{\cal R}(X,Y)Z=(k_1(X\wedge Y-X\wedge R Y-R X\wedge Y) +(k_1+k_2) R X\wedge R Y)Z\vspace{1ex}\\
 \hspace*{13ex}+A_{\alpha(Y,Z)}X -A_{\alpha(X,Z)}Y,\end{array}\ee
 \be\label{eq:codazzi}
 (\nabla^\perp_X \alpha)(Y,Z)-(\nabla^\perp_Y \alpha)(X,Z)=\<\Phi X,Z\>SY-\<\Phi Y,Z\>SX\ee
 and
 \be\label{eq:ricci}
 {\cal R}^\perp(X,Y)\eta=  \alpha(X,A_\eta Y) -\alpha(A_\eta X,Y)+(k_1+k_2)(SX\wedge SY)\eta.
 \ee
 where $\Phi=k_1I-(k_1+k_2)R$. The Codazzi equation (\ref{eq:codazzi}) can also be written as 
 \be\label{eq:codazzi2}(\nabla_Y A)(X,\xi)-(\nabla_X A)(Y,\xi)=\<SX, \xi\>\Phi Y-\<SY, \xi\>\Phi X.\ee

We use the fact that $\Q_{k}^N$, $k\neq 0$, admits a canonical isometric embedding  in  $\R_{\sigma(k)}^{N+1}$ as (a connected component of, if $k<0$)
$$\Q_{k}^N=\{X\in \R_{\sigma(k)}^{N+1}\,:\,\<X,X\>={1}/{k}\}.$$
Here, for $k\in \R$ we set  $\sigma(k)=1$ if $k<0$ and  $\sigma(k)=0$ otherwise, and as a subscript of an Euclidean space it means the index of the corresponding flat metric.
 Thus, 
$\Q_{k_1}^{n_1}\times \Q_{k_2}^{n_2}$ admits
a canonical isometric embedding
\be\label{eq:j}
h\colon\,\Q_{k_1}^{n_1}\times \Q_{k_2}^{n_2}\to\R_{\sigma(k_1)}^{N_1}\times \R_{\sigma(k_2)}^{N_2}=\R_{\mu}^{N_1+N_2},
\ee
with $\mu=\sigma(k_1)+\sigma(k_2)$, $N_i=n_i+1$ if $k_i\neq 0$ and  $N_i=n_i$   otherwise, in which case $\Q_{k_i}^{n_i}$  stands for $\R^{n_i}$. 

Denote by $\tilde\pi_i\colon\,\R_{\mu}^{N_1+N_2}\to \R_{\sigma(k_i)}^{N_i}$  the canonical projection, $i=1,2$. Then, the normal space of $h$ at each point $z\in \Q_{k_1}^{n_1}\times \Q_{k_2}^{n_2}$ is spanned by $k_1\tilde\pi_1(h(z))$ and $k_2\tilde\pi_2(h(z))$, and the second fundamental form of
$h$ is given by
\be\label{eq:sffj} \alpha_h(X,Y)=-k_1\<\pi_1 X, Y\>\tilde\pi_1\circ h-k_2\<\pi_2X, Y\>\tilde\pi_2\circ h.\ee

\vspace{1ex}

Given an isometric immersion $f\colon\,M\to \Q_{k_1}^{n_1}\times \Q_{k_2}^{n_2}$,  set $F=h\circ f$. If $k_i\neq 0$, write $k_i=\epsilon_i/r_i^2$, where $\epsilon_i$ is either $1$ or $-1$, according as $k_i>0$ or $k_i<0$, respectively.
 If $k_1\neq 0$, then the unit  vector field  $\nu_1=\nu_1^{F}=\frac{1}{r_1}\tilde\pi_1\circ F$ is  normal  to $F$, and 
$$\tilde{\nabla}_X\nu_1=\frac{1}{r_1}\tilde\pi_1F_*X=\frac{1}{r_1}(F_*X-h_*LX),$$
where $\tilde \n$ stands for the derivative in $\R_{\mu}^{N_1+N_2}$.
 Hence
 \be\label{eq:nu}^F\nabla_X^\perp \nu_1=-\frac{1}{r_1}h_*SX\,\,\,\,\,\mbox{and}\,\,\,\,\,\,A^{F}_{\nu_1} = -\frac{1}{r_1}(I-R).\ee
 If  $k_2\neq 0$, then $\nu_2=\nu_2^{\tilde{f}}=\frac{1}{r_2}\tilde\pi_2\circ F$ is  also a unit normal vector field   to $F$ such that
$$\tilde{\nabla}_X\nu_2=\frac{1}{r_2}\tilde\pi_2F_*X=\frac{1}{r_2}h_*LX.$$
 Thus
 \be\label{eq:nu'}^F\nabla_X^\perp \nu_2=\frac{1}{r_2}h_*SX\,\,\,\,\,\mbox{and}\,\,\,\,\,\,A^{F}_{\nu_2} = -\frac{1}{r_2}R.\ee
 
 Set \be\label{eq:vartheta1}\vartheta=-\frac{\e_1}{r_1}\nu_1+\frac{\e_2}{r_2}\nu_2.\ee Then $\<\vartheta, \vartheta\>=\frac{\e_1}{r^2_1}+\frac{\e_2}{r^2_2}=k_1+k_2$ and  $\<\vartheta, F\>=0$. Moreover, 
\be\label{eq:vartheta2} A^F_{\vartheta}=\Phi\,\,\,\,\,\mbox{and}\,\,\,\,\,^F\nabla^\perp_X\vartheta=(k_1+k_2)h_*SX.\ee

For later use we prove the following fact.

\begin{proposition}\po\label{prop:comp} Let $f\colon\, M\to N$ and $g\colon\, N\to \Q_{k_1}^{n_1}\times \Q_{k_2}^{n_2}$  be  isometric immersions. Then the tensors $R^g$, $S^g$ and $T^g$ of $g$ and $R^F$, $S^F$ and $T^F$ of $F=g\circ f$ are related by
\be\label{eq:rel1}\<R^FX, Y\>=\<R^gf_*X, f_*Y\>,\,\,\,\,\<S^FX, g_*\xi\>=\<R^gf_*X, \xi\>,\,\,\,\,\,\<S^FX, \zeta\>=\<S^gf_*X, \zeta\>,\ee
\be\label{eq:rel2}\<T^Fg_*\eta, g_*\xi\>=\<R^g\eta, \xi\>, \,\,\,\,\,\<T^Fg_*\eta,\zeta\>=\<S^g\eta, \zeta\>,\,\,\,\,\<T^F\zeta, \beta\>=\<T^g\zeta, \beta\>\ee
for all $X, Y\in TM$, $\xi, \eta\in N_fM$ and $\zeta, \beta\in N_gN$.
\end{proposition}
\proof Notice that $N_FM=g_*N_fM\oplus N_gN$. Given $X\in TM$ we have 
$$g_*f_*R^FX+S^FX=F_*R^FX+S^FX=\pi_2F_*X=\pi_2g_*f_*X=g_*R^gf_*X+S^gf_*X.$$
Taking the inner product of both sides of the preceding equation with $F_*Y=g_*f_*Y$, $g_*\xi$ and $\zeta$, respectively, gives the three equations in 
(\ref{eq:rel1}). For $\eta\in N_fM$ we have
$$g_*R^g\eta+S^g\eta=\pi_2g_*\eta=F_*(S^F)^tg_*\eta+T^Fg_*\eta=g_*f_*(S^F)^tg_*\eta+T^Fg_*\eta.$$
Taking the inner product of the above equation with $F_*Y=g_*f_*Y$ gives again the second equation in  (\ref{eq:rel1}). Taking the inner product with $g_*\xi$ and $\zeta$ gives the first two relations in (\ref{eq:rel2}). To prove the last equation in (\ref{eq:rel2}),  for any $\zeta\in N_gN$ write
$$g_*(S^g)^t\zeta+T^g\zeta=\pi_2\zeta=F_*(S^F)^t\zeta+T^F\zeta=g_*f_*(S^F)^t\zeta+T^F\zeta$$
 and  take the inner product of both sides with $\beta$. Notice that taking the inner product with $F_*Y$ and $g_*\eta$ gives again the third equation in (\ref{eq:rel1}) and the second one in (\ref{eq:rel2}).\vspace{1ex}\qed

We also state as a lemma the following observation.

\begin{lemma}\po\label{le:comp} Let $f\colon\, M\to \Q_{k_1}^{m_1}\times \Q_{k_2}^{m_2}$ be an isometric immersion and let  $j\colon\, \Q_{k_1}^{m_1}\times \Q_{k_2}^{m_2}\to \Q_{k_1}^{n_1}\times \Q_{k_2}^{n_2}$  be  a totally geodesic inclusion. Set $F=j\circ f$. Then the tensors $R^f$, $S^f$ and $T^f$ of $f$ and $R^F$, $S^F$ and $T^F$ of $F$ are related by
\be\label{eq:rela}R^F=R^f, \,\,\,\,S^F=j_*S^f\,\,\,\mbox{and}\,\,\,T^Fj_*=j_*T^f\ee
\end{lemma}

We conclude this section by  listing   well known formulae for the second fundamental form and normal connection of a composition of isometric immersions. We omit the proofs, which follow by  a straightforward application of the Gauss and 
Weingarten formulae.

\begin{proposition}\po\label{prop:comp2} Let $f\colon\, M\to N$ and $g\colon\, N\to P$ be  isometric immersions. Set $F=g\circ f$. Then 
$N_FM=g_*N_fM\oplus N_gN$ and the second fundamental forms and normal connections of $f$, $g$ and $F$ are related by
\be\label{eq:relA} \alpha_F(X,Y)=g_*\alpha_f(X,Y)+\alpha_g(f_*X,f_*Y),\ee
\be\label{eq:relB} ^F\nabla_X^\perp g_*\xi=g_*\,^f\nabla_X^\perp \xi +\alpha_g(f_*X,\xi),\,\,\,\,\,\,\,(^F\nabla_X^\perp \zeta)_{N_gN}=\,^g\nabla_{f_*X}^\perp \zeta\ee
and
\be\label{eq:relC}\<^F\nabla_X^\perp \zeta,g_*\xi\>=-\<A^g_\zeta f_*X, \xi\>\ee
for all $X\in TM, \xi\in N_fM$ and $\zeta\in N_gN$.
\end{proposition}

\section{Products of isometric immersions}

We start this section by characterizing  isometric immersions into a slice of $\Q_{k_1}^{n_1}\times \Q_{k_2}^{n_2}$ as those for which either  $R$ or $I-R$ vanishes identically.
 
 \begin{proposition}\po\label{prop:R=0} Let $f\colon\,M^m\to \Q_{k_1}^{n_1}\times \Q_{k_2}^{n_2}$  be an isometric immersion. Then 
 $f(M^m)\subset \Q_{k_1}^{n_1}\times \{x_2\}$ for some  $\{x_2\}\in \Q_{k_2}^{n_2}$ if and only if $R=0$.\end{proposition}
 \proof We have $f(M^m)\subset \Q_{k_1}^{n_1}\times \{x_2\}$ for some  $\{x_2\}\in \Q_{k_2}^{n_2}$ if and only $\pi_2\circ f_*=0$, which is equivalent to $R=0$.\vspace{1ex}\qed
 
 Next we show that
  products
 $$f=f_1\times f_2\colon\, M^m=M_1^{m_1}\times M_2^{m_2}\to \Q_{k_1}^{n_1}\times \Q_{k_2}^{n_2}$$
 of isometric immersions  are characterized by 
 having vanishing tensor $S$.

 \begin{lemma}\po\label{le:kerS} Let $f\colon\,M^m\to \Q_{k_1}^{n_1}\times \Q_{k_2}^{n_2}$  be an isometric immersion. Then $\ker S$ splits orthogonally pointwise as 
 $$\ker S=\ker R\oplus \ker(I-R).$$
Moreover, the following holds:
 \begin{itemize}
 \item[$(i)$] If $\ker S$ has constant dimension, then it is a smooth subbundle of $TM$ and so are $\ker R$ and $\ker(I-R)$;
 \item[$(ii)$] If $S=0$ then $\ker R$ and $\ker(I-R)$ are parallel subbundles of $TM$.
 \end{itemize}
 \end{lemma}
 \proof The first assertion follows from the first equation in (\ref{eq:pi1}). Assume that $\ker S=(\Ima S^t)^\perp$ has constant dimension $k$ on $M^m$. Given $x\in M^m$, let $\xi_1, \ldots, \xi_{m-k}$ be normal vectors at $x$ such that $\{S^t\xi_i\}_{1\leq i\leq m-k}$ is linearly independent. Extend $\xi_1, \ldots, \xi_{m-k}$ to smooth normal vector fields on a neighborhood of $x$. Then $\{S^t\xi_i\}_{1\leq i\leq m-k}$ is still linearly independent on a (possibly smaller) neighborhood $V$ of $x$, hence it spans the image of $S^t$ on $V$. It follows that $\Ima S^t$, and hence $\ker S$, is a smooth distribution. Moreover, using the lower semicontinuity of the rank of both $R$ and $I-R$ we easily obtain  that   $A_\ell=\{x\in M\,:\,\ker R(x)=\ell\}$ is an open subset of $M^m$  for any $0\leq \ell\leq m$. Hence $M^m=A_\ell$ for some $\ell\in \{0,\ldots, m\}$, that is, $\ker R$ has constant dimension $\ell$ on $M^m$.  Arguing as before we conclude that $\ker R=\Ima (I-R)$ and $\ker (I-R)=\Ima R$ are also smooth. Finally, from (\ref{eq:pid1}) we obtain that $\nabla R=0$ if $S=0$, and  $(ii)$ follows.\qed

\begin{proposition}\po\label{prop:S=0} Let $f\colon\,M^m\to \Q_{k_1}^{n_1}\times \Q_{k_2}^{n_2}$  be an isometric immersion. Then 
 $M^m$ is locally a  Riemannian product $M^m=M_1^{m_1}\times M_2^{m_2}$ and $f=f_1\times f_2$, where $f_i\colon\,M^{m_i}_i\to \Q_{k_i}^{n_i}$,  $1\leq i\leq 2$, is an isometric immersion,  if and only if  $S=0$ and neither $R=0$ nor $R=I$. The ``if"  part holds globally if $M^n$ is complete and simply connected.\end{proposition}
\proof   Assume that $M^m$ is locally a  Riemannian product $M^m=M_1^{m_1}\times M_2^{m_2}$ and $f=f_1\times f_2$, where $f_i\colon\,M^{m_i}_i\to \Q_{k_i}^{n_i}$,  $1\leq i\leq 2$, is an isometric immersion. Then $TM$ splits orthogonally as $TM=TM_1\oplus TM_2$, with  $f_*TM_i={f_i}_*TM_i\subset T\Q_{k_i}^{n_i}$, $1\leq i\leq 2$. Hence
$R=0$ on $TM_1$ and $R=I$ on $TM_2$. Thus $S^tS=R(I-R)=0$, and therefore $S=0$. 

Conversely, suppose that $S=0$ and that neither $R=0$ nor $R=I$. Then  $\ker R$ and $\ker(I-R)$ are nontrivial parallel subbundles of $TM$, and $TM$ splits orthogonally as $TM=\ker R \oplus\ker(I-R)$ by  Lemma \ref{le:kerS}. That $M^m$ splits locally as a Riemannian product $M^m=M_1\times M_2$  then follows from the local version of de Rham theorem.

Since $S=0$, formula (\ref{eq:pid2}) becomes
$$\alpha_f(RX,Y)=T\alpha_f(X,Y)$$
for all $X, Y\in TM$.
The right-hand-side is symmetric with respect to $X$ and $Y$, so the same holds for the left-hand-side, i.e.,
$$ \alpha_f(RX,Y)=\alpha_f(X,RY)$$
for all $X, Y\in TM$. It follows that $\alpha_f(X,Y)=0$ if $X\in \ker R$ and $Y\in \ker(I-R)$. 

Set $F=h\circ f$, where $h$ is the inclusion of $\Q_{k_1}^{n_1}\times \Q_{k_2}^{n_2}$ into $\R_\mu^{N_1+N_2}=\R_{\sigma(k_1)}^{N_1}\times \R_{\sigma(k_2)}^{N_2}$ as in (\ref{eq:j}). Then the second formulas in (\ref{eq:nu}) and 
(\ref{eq:nu'}) show that we also have  $\alpha_{F}(X,Y)=0$ if $X\in \ker R$ and $Y\in \ker(I-R)$. Define
 $$V_1=\spa\{F_*(x)X:x\in M^m,\, X\in \ker R(x)\}$$ and 
$$V_2=\spa\{F_*(x)X:x\in M^m,\, X\in \ker(I-R(x))\}.$$ Since $\tilde\pi_2(F_*\ker R)=\{0\}=\tilde\pi_1(F_*\ker(I-R))$, we have   $V_1\subset \R_{\sigma(k_1)}^{N_1}$ and $V_2\subset \R_{\sigma(k_2)}^{N_2}$. As in the proof of Moore's Lemma \cite{mo}, it follows that there exist isometric immersions $F_1\colon\,M_1\to \R_{\sigma(k_1)}^{N_1}$ and $F_2\colon\,M_2\to \R_{\sigma(k_2)}^{N_2}$ such that $F(x,y)=(F_1(x), F_2(y))$. Since $$F(M_1\times M_2)=F_1(M_1)\times F_2(M_2)\subset \Q_{k_1}^{n_1}\times \Q_{k_2}^{n_2},$$ we have that $F_i=j_i\circ f_i$ for some isometric immersions $f_i\colon\, M_i\to \Q_{k_i}^{n_i}$, $1\leq i\leq 2$,   where $j_i$ is the  inclusion of $\Q_{k_i}^{n_i}$ into $\R_{\sigma(k_i)}^{N_i}$.

  The global assertion follows as before from  the global version of de Rham Theorem.\qed
  
  \section{A reduction of codimension theorem}
  
 We say that an isometric immersion $f\colon\, M^m\to \Q_{k_1}^{n_1}\times \Q_{k_2}^{n_2}$ \emph{reduces codimension on the left by $\ell$} if there exists a totally geodesic inclusion $j_1\colon\, \Q_{k_1}^{m_1}\to \Q_{k_1}^{n_1}$ with $n_1-m_1=\ell$ and an isometric immersion $\bar f\colon\, M^m\to \Q_{k_1}^{m_1}\times \Q_{k_2}^{n_2}$ such that $f=(j_1\times id)\circ \bar f$. Similarly we define what it means by $f$ reducing codimension \emph{on the right}. In this section we give necessary and sufficient conditions for an isometric immersion $f\colon\, M^m\to \Q_{k_1}^{n_1}\times \Q_{k_2}^{n_2}$ to reduce codimension on the left or on the right. We start with the following observation.
 
 \begin{lemma}\po\label{le:stmperp} Let $ f\colon\, M^m\to \Q_{k_1}^{n_1}\times \Q_{k_2}^{n_2}$ be an isometric immersion. Then 
 $$S(TM)^\perp=U\oplus V,\,\,\,\,\mbox{where}\,\,\,U=\ker T\,\,\,\mbox{and}\,\,\,\,V=\ker (I-T).$$
 \end{lemma}
 \proof It follows from the third equation in (\ref{eq:pi1}) that  $\ker T(I-T)=\ker S^t=S(TM)^\perp$, hence the restriction of  $T$ to $S(TM)^\perp$ is the orthogonal projection onto $V$.\vspace{1ex}\qed

The following result and its corollary are the analogues for isometric immersions into products of space forms of the well known criterion for reduction of codimension of isometric immersions into space forms (see, e.g., Proposition $4.1$ and Corollary $4.2$ in Chapter $4$ of \cite{daa}).
We restrict ourselves to stating the results for reduction of codimension \emph{on the left}, for the case of reduction of codimension on the right is completely similar, just by replacing the vector subbundle $U$ by $V$ in what follows. 

Given an isometric immersion $f\colon\, M\to N$ between Riemannian manifolds, its \emph{first normal space} at $x\in M$ is the subspace $N_1(x)$ of $N_fM(x)$ spanned by the image of its second fundamental form $\alpha_f(x)$. 

\begin{theorem}\po\label{thm:redcod1} Let $f\colon\,M^m\to \Q_{k_1}^{n_1}\times \Q_{k_2}^{n_2}$ be an isometric immersion. Then the following assertions are equivalent:
\begin{itemize}
\item[$(i)$] $f$ reduces codimension on the left by $\ell$;
\item[$(ii)$] There exists a subbundle $L$ of rank $\ell$ of $N_fM$ such that $L$ is parallel in the normal connection and $L\subset U\cap N_1^\perp$.
\end{itemize}
\end{theorem} 
\proof  Assume that there exists a totally geodesic inclusion $j_1\colon\, \Q_{k_1}^{m_1}\to \Q_{k_1}^{n_1}$ with $\ell=n_1-m_1$ and an isometric immersion $\bar f\colon\, M^m\to \Q_{k_1}^{m_1}\times \Q_{k_2}^{n_2}$ such that $f=j\circ \bar f$, where $j=j_1\times id$. 
Set $$L=N_{j_1}\Q_{k_1}^{m_1}=N_j(\Q_{k_1}^{m_1}\times \Q_{k_2}^{n_2})\subset N_{f}M.$$ Since $j$ is totally geodesic, it follows  from  (\ref{eq:relA}) (with $\bar f$, $j$ and $f$ playing the roles of $f$, $g$ and $F$, respectively) that $N_1\subset L^\perp$, hence $L\subset N_1^\perp$. Also because $j$ is totally geodesic, we obtain from  (\ref{eq:relC}) that $L$ is parallel in the normal connection. Finally, since $\pi_2$ vanishes on $L=N_j\Q_{k_1}^{m_1}$, it follows that $L\subset U$.

Now we prove the converse. Set $F=h\circ f$, where $h$ is as in (\ref{eq:j}). For any $\xi\in N_fM$, using (\ref{eq:sffj})  we obtain
\be\label{eq:der}
\tilde \nabla_Xh_*\xi=h_*\nabla^\perp_X\xi-F_*A_\xi X+\alpha_h(f_*X, \xi)=h_*\nabla^\perp_X\xi-F_*A_\xi X-\<SX,\xi\>\vartheta,
\ee
with $\vartheta$ as in (\ref{eq:vartheta1}). Given $\xi\in L$, using that $L$ is parallel in the normal connection and that $L\subset U\cap N_1^\perp$ it follows  from   (\ref{eq:der}) that 
\be\label{eq:W2} \tilde \nabla_Xh_*\xi=h_*\nabla^\perp_X\xi\in h_*L.\ee
This shows that  $W:=h_*L$ is a  constant subspace in $\R_{\mu}^{N_1+N_2}$. Moreover, we have that $\tilde \pi_2h_*\xi=h_*\pi_2\xi=h_*T\xi=0$, because $\pi_2|_{S(TM)^\perp}=T$ and $\xi\in U=\ker T$. Hence  $\tilde \pi_2|_W=0$.\vspace{1ex}

\noindent  {\bf Claim:} $\tilde{\pi}_1(F(M^m))\subset \Q_{k_1}^{n_1}\cap W^\perp$. \vspace{1ex}

Given $\xi \in L(x)$, $x\in M^n$ and $X\in T_xM^n$, we have using (\ref{eq:W2}) that
$$X\<\tilde{\pi}_1\circ F, h_*\xi\>=\<\tilde{\pi}_1 F_*X,h_*\xi\>=-\<SX,\xi\>=0,$$
since $L\subset U\subset S(TM)^\perp$. 
Hence $\tilde{\pi}_1(F(M^n))\subset \tilde{\pi}_1(F(x_0))+W^\perp$ for any fixed $x_0\in M^n$. But $\tilde{\pi}_1(F(x_0))\in W(x_0)^\perp=W^\perp$, hence $\tilde{\pi}_1(F(M^n))\subset \Q_{k_1}^{n_1}\cap W^\perp$ as claimed. \vspace{1ex}\qed

\begin{corollary}\po\label{cor:redcod} Let $f\colon\,M^m\to \Q_{k_1}^{n_1}\times \Q_{k_2}^{n_2}$ be an isometric immersion. Assume that $U\cap N_1^\perp$ is a vector subbundle of $N_fM$ with rank $\ell$  satisfying
\be\label{eq:condred}\nabla^\perp (U\cap N_1^\perp)\subset N_1^\perp.\ee
 Then $f$ reduces codimension on the left by $\ell$.
\end{corollary} 
\proof By Theorem \ref{thm:redcod1}, it suffices to prove that 
\be\label{eq:un1}\nabla^\perp (U\cap N_1^\perp)\subset U.\ee
 But this follows from  (\ref{eq:pid4}),  for it implies for all $\xi\in U\cap N_1^\perp$ and $X\in TM$ that
$$T\nabla_X^\perp \xi=\nabla_X^\perp T\xi=0.\vspace{1ex}\qed$$

In case $U\cap N_1^\perp$ is a vector subbundle of the normal bundle of an isometric immersion $f\colon\,M^n\to \Q_{k_1}^{n_1}\times \Q_{k_2}^{n_2}$, the next result gives necessary and sufficient conditions for (\ref{eq:condred}) to hold in terms of its normal curvature tensor and mean curvature vector field. It is the version for submanifolds of products of space forms of a theorem by Dajczer for submanifolds of space forms (see \cite{da} or Theorem $4.4$ in \cite{daa}). 

\begin{theorem}\po \label{thm:dajczer} Let $f\colon\,M^m\to \Q_{k_1}^{n_1}\times \Q_{k_2}^{n_2}$ be an isometric immersion. Assume that $U\cap N_1^\perp$ is a vector subbundle of $N_fM$.
Then $\nabla^\perp (U\cap N_1^\perp)\subset N_1^\perp$
 if and only if the following conditions hold:
	\begin{itemize}
	\item[$(i)$] $\nabla^\perp \mathcal R^\perp|_{U\cap {N_1}^\perp}=0$,
	\item[$(ii)$] $\nabla^\perp (U\cap N_1^\perp)\subset \{H\}^\perp$.
	\end{itemize}
\end{theorem}
\proof Assume that $\nabla^\perp (U\cap N_1^\perp)\subset N_1^\perp$. Then $(ii)$ is clear. 
 By the Ricci equation (\ref{eq:ricci}), for any $\xi \in U\cap {N_1}^\perp\subset S(TM)^\perp\cap N_1^\perp$ we have
$$\mathcal{R}^\perp(X,Y)\xi = \alpha(X, A_\xi Y) - \alpha(A_\xi X,Y) + (k_1 + k_2) (S X \wedge S Y) \xi =0.
$$
Using that $\nabla^\perp_Z \xi \in U\cap N_1^\perp$ by (\ref{eq:un1}) and the assumption,  we obtain
$$	\begin{array}{l}\left(\nabla^\perp_Z \mathcal{R}^\perp \right)(X,Y,\xi) =\nabla^\perp_Z \mathcal{R}^\perp(X,Y) \xi - \mathcal{R}^\perp(\nabla_Z X, Y) \xi - \mathcal{R}^\perp(X, \nabla_Z Y) \xi -\vspace{1ex}\\\hspace*{21ex} -\mathcal{R}^\perp(X, Y) \nabla^\perp_Z \xi = 0.
	\end{array}$$
	Hence $\nabla^\perp \mathcal R^\perp|_{U\cap {N_1}^\perp}=0$.
	
We now prove the converse.  Let  $\xi \in  U\cap N_1^\perp$. We must prove that $\nabla^\perp_Z\xi\in N_1^\perp$ for all $Z\in TM$. We obtain from  $(i)$  that  $\mathcal{R}^\perp (X,Y) \nabla^\perp_Z\xi = 0$. Since $\nabla^\perp_Z \xi \in U\subset S(TM)^\perp$ by (\ref{eq:un1}), the Ricci equation (\ref{eq:ricci}) yields

$$ \left[A_{\nabla^\perp_Z \xi}, A_{\nabla^\perp_W \xi} \right] = 0\,\,\,\mbox{for all}\,\,\, W,Z \in TM.$$

	Hence, at any $x \in M$ there exists an  orthonormal basis $\left\{E_1(x), \cdots, E_m(x)\right\}$ of $T_x M$ that  diagonalizes simultaneously the family of operators $\left\{ A_{\nabla^\perp_X \xi}\,:\,X \in T_x M\right\}$.  Then
	$$\<\alpha({E_i},{E_j}),\nabla^\perp_{E_k} \xi\> = \<A_{\nabla^\perp_{E_k}\xi} E_i,{E_j}\> = 0,\,\,\,\mbox{if}\,\,\,i\neq j.
	$$
Using  that $\xi \in U\cap N_1^\perp\subset S(TM)^\perp\cap N_1^\perp$, we obtain from the Codazzi equation (\ref{eq:codazzi2})  that 
$$	A_{\nabla^\perp_X \xi} Y = A_{\nabla^\perp_Y \xi} X\,\,\,\mbox{for all} \,\,\,X, Y \in TM,$$
which implies that $A_{\nabla^\perp_{E_i} \xi} E_j =0$ for $i\neq j$. 
Hence 
	$$\<\alpha({E_j},{E_j}),\nabla^\perp_{E_k} \xi\> =\<A_{\nabla^\perp_{E_k} \xi} E_j, E_j\>=0
	$$
	if $k\neq j$. Finally, using condition $(ii)$ and the above we obtain
	$$\<\alpha({E_j},{E_j}),\nabla^\perp_{E_j} \xi\> =n\<H,\nabla^\perp_{E_j} \xi\>
	=0.$$
	Hence $\nabla^\perp_Z\xi\in N_1^\perp$ for all $Z\in TM$. 
\qed

\begin{remark}\po {\em Given an isometric immersion  $f\colon\,M^m\to \Q_k^n\times\R$,  
 let   
 $\frac{\partial}{\partial t}$ be a unit
vector field tangent to the second factor. 
Then, a tangent vector field  $Z$ on $M^m$ and a normal vector field  $\eta$ along $f$ are defined by
$$\frac{\partial}{\partial t}={f}_* Z+ \eta.$$
The tensors $R$, $S$ and $T$ associated to $f$ are given by
$$RX=\<X,Z\>Z, \,\,\,\,\,\,SX=\<X,Z\>\eta\,\,\,\,\,\mbox{and}\,\,\,\,\,T\xi=\<\xi, \eta\>\eta.$$
Then $U=\ker T=\{\eta\}^\perp$, hence $U\cap N_1^\perp =(N_1+ \spa\{\eta\})^\perp$. Thus,  condition (\ref{eq:condred}) is equivalent to
$$\nabla^\perp N_1\subset N_1+\spa\{\eta\}.$$
Therefore, in this case Corollary \ref{cor:redcod} and Theorem \ref{thm:dajczer} reduce to Lemma $6$ and Theorem~$7$ of \cite{mt}, respectively.}
\end{remark}

\section{Weighted sums of isometric immersions}
 
 Given  $a, b\in \R^*$ with $a^2+b^2=1$, set $\tilde{k}_1=a^2 k_1$ and $\tilde k_2=b^2k_2$. Let $f_i\colon\, M^m\to \Q_{\tilde k_i}^{n_i}\subset \R_{\sigma(k_i)}^{N_i}$ be isometric immersions, $1\leq i\leq 2$. Then  $f=(af_1,bf_2)\colon\,M^m\to \R_{\mu}^{N_1+N_2}$ is an isometric immersion 
 that takes values in $\Q_{k_1}^{n_1}\times  \Q_{k_2}^{n_2}$. We call $f$ the \emph{weighted sum of $f_1$  and $f_2$} with weights $a$ and $b$.  The normal space  of $f$ (in $\R_{\mu}^{N_1+N_2}$) is the orthogonal sum
 $$N_fM=N_{f_1}M\oplus N_{f_2}M\oplus W,$$
 where $W=\spa\{-b{f_1}_*X+a{f_2}_*X\,:\, X\in TM\}$.
 We have 
 $$(\pi_2\circ f_*)X=b{f_2}_*X=f_*b^2X+ab(-b{f_1}_*X+a{f_2}_*X),$$ 
 with $-b{f_1}_*X+a{f_2}_*X\in N_fM$, hence $R=b^2I$ and $SX=ab(-b{f_1}_*X+a{f_2}_*X)$  for any $X\in TM$. In particular,
 we have $S(TM)= W$. 
 
 \begin{example}\label{ex:tgeod} {\em Let $k_1, k_2\in \R$ with $k_1k_2>0$ and  set $$a:=\sqrt{\frac{k_2}{k_1+k_2}},\,\,  b:=\sqrt{\frac{k_1}{k_1+k_2}}\,\,\,\mbox{and}\,\,\,\e=\sigma(k_1)=\sigma(k_2).$$ Given  $T_i\in O_{\e}(n+1)$, $1\leq i\leq 2$, define  $G\colon\, \R_\e^{n+1}\to \R_{2\e}^{2n+2}=\R_\e^{n+1}\oplus \R_\e^{n+1}$ by $$G(x)=(aT_1(x), bT_2(x)).$$ Then $G(\Q_k^n)\subset \Q_{k_1}^{n}\times  \Q_{k_2}^{n}$, $k=k_1k_2/(k_1+k_2)$, thus $G|_{\Q_k^n}=h\circ g$ for some  isometric immersion $g\colon\,\Q_k^n\to \Q_{k_1}^{n}\times  \Q_{k_2}^{n}$, where $h$ is as in (\ref{eq:j}). Since  $G(\Q_k^n)=V\cap \Q_k^{2n+1}$, where $V=G(\R_\e^{n+1})$, it follows that $\bar h\circ g$ is totally geodesic, where $\bar h\colon\,\Q_{k_1}^{n}\times  \Q_{k_2}^{n}\to \Q_k^{2n+1}$ is the inclusion. Hence $g$ is totally geodesic. Moreover, $R^g=b^2I$ by the preceding discussion.}
 \end{example}
 
 The following result shows that weighted sums of isometric immersions are characterized by the fact that the tensor $R$ is a multiple of the identity tensor.
 
 \begin{proposition}\label{prop:rmid}  Let $f\colon\, M^m\to \Q_{k_1}^{n_1}\times \Q_{k_2}^{n_2}$  be an isometric immersion.  Then the following assertions are equivalent:
 \begin{itemize}
 \item[$(i)$]  there exist  isometric immersions $f_i\colon\, M^m\to \Q_{\tilde k_i}^{n_i}\subset \R_{\sigma(k_i)}^{N_i}$, $1\leq i\leq 2$, with   $\tilde{k}_1=k_1\cos^2\theta$ and $\tilde k_2=k_2\sin^2\theta$ for some $\theta\in (0, \pi/2)$,  such that 
 $f=(\cos\theta f_1,\sin \theta f_2)$;
 \item[$(ii)$] $R=\sin^2\theta I$ for some $\theta\in (0, \pi/2)$.
 \end{itemize}
 \end{proposition}
 \proof We have that $R=\sin^2\theta I$ for some $\theta\in (0, \pi/2)$ if and only if the tensor $L=\pi_2\circ f_*$ satisfies
 $$L^tL=R=\sin^2\theta I.$$
 This is equivalent to  $\pi_2\circ f_*$ being a  similarity of ratio $\sin \theta$. In turn, this holds if and only if  there exist  isometric immersions $f_i\colon\, M^m\to \Q_{\tilde k_i}^{n_i}\subset \R^{n_i+1}$, $1\leq i\leq 2$, with $\tilde{k}_1=k_1\cos^2\theta$ and $\tilde k_2= k_2\sin^2 \theta$, such that $\pi_1\circ f=\cos \theta f_1$ and $\pi_2\circ f=\sin \theta f_2$. \qed
 
 \section{A further theorem on reduction of codimension}
 
  We derive in this section necessary and sufficient conditions for the image of an isometric immersion 
  $f\colon\, M^m\to \Q_{k_1}^{n_1}\times \Q_{k_2}^{n_2}$ to be contained in the totally geodesic submanifold of Example \ref{ex:tgeod}.
  This will be used in the proof of Theorem \ref{thm:parallel} in the next section. 
  We need some preliminary results.
 
 \begin{lemma}\po\label{le:nablaR} Let $f\colon\, M^m\to \Q_{k_1}^{n_1}\times \Q_{k_2}^{n_2}$  be an isometric immersion. 
 Then $\nabla R=0$ if and only if $S(TM)\subset N_1^\perp$.
 \end{lemma}
 \proof By (\ref{eq:pid1}), we have that $\nabla R=0$ if and only if 
 $$A_{SY}X+S^t\alpha(X,Y)=0$$
 for all $X, Y\in TM$. This is equivalent to
 $$\<\alpha(X, Z), SY\>=-\<\alpha(X,Y), SZ\>$$
 for all $X,Y,Z\in TM$. Hence, $\nabla R=0$ if and only if the trilinear form $\beta$ given by
 $$\beta(X,Y,Z)=\<\alpha(X,Y), SZ\>$$
 is skew-symmetric in the last two variables. Since $\beta$ is symmetric in the first two variables, it follows from the next lemma, called the \emph{Braid Lemma} (see Section 9.5.4.9 of \cite{be}), that this is the case if and only if $\beta$ vanishes.\qed
 
 \begin{lemma} Let $\beta\colon\, V\times V\times V\to W$ be a trilinear map. If $\beta$ is symmetric in the first two variables and skew-symmetric in the last two, then $\beta=0$.
 \end{lemma} 
 \proof For any $X, Y, Z\in V$ we have
 \begin{eqnarray*}\beta(X,Y,Z)&=&-\beta(X,Z,Y)=-\beta(Z,X,Y)=\beta(Z,Y,X)=\beta(Y,Z,X)=-\beta(Y,X,Z)\\
 &=&-\beta(X,Y,Z),
 \end{eqnarray*}
 hence $\beta=0$.\qed
 
 \begin{lemma} \po \label{le:R} Let $f\colon\, M^m\to \Q_{k_1}^{n_1}\times \Q_{k_2}^{n_2}$ be an isometric immersion.  If $L$ is a subbundle of $S(TM)^\perp$ such that $$(T-\lambda I)L\subset L^\perp\,\,\,\mbox{for some}\,\,\,\lambda\in \R,$$ then ($\lambda\in [0,1]$ and ) $\pi_2|_L$ is a similarity of ratio $\lambda$.
 \end{lemma}
 \proof  Since $L\subset S(TM)^\perp$ and  $T|_{S(TM)^\perp}=\pi_2|_{S(TM)^\perp}$, we have for all $\xi, \eta\in L$ that
 $$\<\pi_2\xi, \pi_2\eta\>=\<\pi_2\xi,\eta\>=\<T\xi,\eta\>=\lambda\<\xi,\eta\>.\qed$$

  \begin{lemma}\po\label{le:phi}  Let $f\colon\, M^m\to \Q_{k_1}^{n_1}\times \Q_{k_2}^{n_2}$, $k_1k_2\neq 0$,  be an  isometric immersion with  $\Phi=0$. Then $k_1k_2>0$ and $R=b^2 I$, with $b:=\sqrt{\frac{k_1}{k_1+k_2}}$.
  \end{lemma}
  \proof Since $$0=\Phi=k_1I-(k_1+k_2)R,$$ if $k_1+k_2=0$ then $k_1=0$, in contradiction with the assumption that $k_1k_2\neq 0$. Hence $k_1+k_2\neq 0$ and $R=\lambda I$ with $\lambda=\frac{k_1}{k_1+k_2}\in (0,1)$. In particular,  we have that $\frac{k_1k_2}{(k_1+k_2)^2}=\lambda(1-\lambda)>0$, hence $k_1k_2>0$.\qed

 \begin{theorem}\po\label{thm:redcod2}  Let $f\colon\, M^m\to \Q_{k_1}^{n_1}\times \Q_{k_2}^{n_2}$, $k_1k_2\neq 0$,  be an  isometric immersion. Then the following assertions are equivalent:
\begin{itemize}
 \item[$(i)$] $k_1k_2>0$ and there exist $\ell\geq 0$, an isometric immersion $\bar f\colon\, M^m\to \Q_{k}^{m+\ell}$, with $k=k_1k_2/(k_1+k_2)$, and  a totally geodesic inclusion  $j\colon\,\Q_{k_1}^{m+\ell}\times \Q_{k_2}^{m+\ell}\to \Q_{k_1}^{n_1}\times \Q_{k_2}^{n_2}$ such that $$f=j\circ g\circ \bar f,$$ where $g\colon\,  \Q_{k}^{m+\ell}\to \Q_{k_1}^{m+\ell}\times \Q_{k_2}^{m+\ell}$ is the totally geodesic embedding of Example \ref{ex:tgeod}.
 \item[$(ii)$] $\Phi=0$ and there exists a subbundle $L^\ell$ of $S(TM)^\perp$ such that $N_1\subset L$, $L$ is parallel in the normal connection and $(T-b^2I)L\subset L^\perp$, where $b:=\sqrt{\frac{k_1}{k_1+k_2}}$.
 \end{itemize}
 \end{theorem}
 \proof Let us prove that $(ii)$ implies $(i)$. By Lemma \ref{le:phi}, we have $k_1k_2>0$, hence $\sigma(k_1)=\sigma(k_2):=\epsilon\in \{0, 1\}$. Let $h\colon\, \Q_{k_1}^{n_1}\times \Q_{k_2}^{n_2}\to \R_{2\epsilon}^N$, $N=n_1+n_2+2$,  be the canonical inclusion and set $F=h\circ f$. \vspace{1ex}\\
 {\bf Claim $1$:}  $V:=F_*TM\oplus h_*L\oplus \spa\{F\}$ is a constant subspace of $\R^N_{2\epsilon}$.\vspace{1ex} 
 
 To prove Claim $1$, it suffices to show that the  orthogonal complement $V^\perp$ of $V$ in $\R^N_{2\epsilon}$ is a constant subspace. We have
 $$V^\perp=h_*L^\perp \oplus \spa\{\vartheta\},$$ where $\vartheta$ is as in (\ref{eq:vartheta1}). By the assumption and the first equation in (\ref{eq:vartheta2}) we have $$A^F_{\vartheta}=\Phi=0.$$ 
 Using this and the second equation in (\ref{eq:vartheta2}) we obtain
 $$\tilde \nabla_X\vartheta=-F_*A^F_{\vartheta}X+^F{\!}\nabla^\perp_X\vartheta=(k_1+k_2)h_*SX\in h_*L^\perp \subset V^\perp,$$
 because $S(TM)\subset L^\perp$ by assumption.
 On the other hand, using that $L$ is parallel in the normal connection and that $L^\perp \subset N_1^\perp$ we have
 $$\tilde \nabla_Xh_*\xi=h_*\bar \nabla_X\xi+\alpha_h(f_*X, \xi)=h_*\nabla^\perp_X\xi-\<SX,\xi\>\vartheta\in V^\perp$$
 for all $\xi \in L^\perp$, where $\bar \nabla$ stands for the connection on $\Q_{k_1}^{n_1}\times \Q_{k_2}^{n_2}$. Thus Claim $1$ is proved.
 
 Since  $F_*TM\subset V$, we have that $F(M)\subset V$, because $V$ contains the position vector at any point. \vspace{1ex}\\
 {\bf Claim $2$:} $\tilde \pi_2|_V$ is a similarity of ratio $b$.\vspace{1ex}
 
Since $L\subset S(TM)^\perp$ and $(T-b^2I)L\subset L^\perp$ by assumption, it follows from Lemma \ref{le:R}  that $\tilde \pi_2|_{h_*L}$ is a similarity of ratio $b$. We also have that
 $$\<\tilde\pi_2F, \tilde\pi_2F\>=\frac{1}{k_2}=b^2\frac{k_1+k_2}{k_1k_2}=b^2\<F, F\>$$
 and that $\tilde \pi_2|_{F_*TM}$ is a similarity of ratio $b$, because $R=b^2 I$ by Lemma \ref{le:phi}. The proof of Claim $2$ is completed by noticing that
 $$\<\tilde\pi_2F, \tilde\pi_2F_*X\> =\<\tilde\pi_2F, F_*X\> = 0,$$
 $$\<\tilde\pi_2F,\tilde\pi_2 h_*\xi\> =\<\tilde\pi_2F, h_*\xi\> = 0$$
  and 
 $$\<\tilde\pi_2F_*X, \tilde\pi_2 h_*\xi\> = \<SX, \xi\> = 0$$
 for any $\xi\in L$.

	Let $\ell$ be the rank of $L$, set $a=\sqrt{\frac{k_2}{k_1+k_2}}$ and denote by $\R^{m+\ell+1}$ the image of both $\tilde \pi_1|_V$ and $\tilde \pi_2|_V$ in $\R_\epsilon^{n_1+1}$ and $\R_\epsilon^{n_2+1}$, respectively. By Claim $2$, we have that $T_1=a^{-1}\tilde \pi_1|_V$ and $T_2=b^{-1}\tilde \pi_2|_V$ are linear isometries onto $\R_\epsilon^{m+\ell+1}$, and 
	$$V=\{G(X):=(aT_1(X), bT_2(X))\,:\,X\in \R_\epsilon^{m+\ell+1}\}.$$
	Let $g=G|_{\Q_{k}^{m+\ell}}\colon\, \Q_{k}^{m+\ell}\to \Q_{k_1}^{m+\ell}\times \Q_{k_2}^{m+\ell}$ be the totally geodesic embedding  in Example~\ref{ex:tgeod} with $g(\Q_{k}^{m+\ell})=V \cap \Q_{k}^{2m+2\ell+1}$.

	Since $f(M) \subset g(\Q_{k}^{m+\ell})\subset \Q_{k_1}^{m+\ell}\times \Q_{k_2}^{m+\ell}\subset \Q_{k_1}^{n_1}\times \Q_{k_2}^{n_2}$, it follows that $f = j \circ \tilde{f}$, where $j\colon\,\Q_{k_1}^{m+\ell}\times \Q_{k_2}^{m+\ell}\to \Q_{k_1}^{n_1}\times \Q_{k_2}^{n_2}$ is a totally geodesic inclusion and  $\tilde{f}\colon\, M \to \Q_{k_1}^{m+\ell}\times \Q_{k_2}^{m+\ell}$ is an isometric immersion with  $\tilde{f}(M) \subset g(\Q_{k}^{m+\ell})$. Therefore, $\tilde f=g\circ \bar f$ for some isometric immersion 
	$\bar f\colon\, M^n\to \Q_k^{m+\ell}$. 
	
	For the converse, define $L=j_*g_*N_{\bar f}M$. From formula (\ref{eq:relA}), with $\bar f$, $j\circ g$ and $f$ playing the roles of $f$, $g$ and $F$ in that formula, respectively, it follows that $N_1^f\subset L$, because $j\circ g$ is totally geodesic. That $L$ is parallel in the normal connection follows from the first formula in (\ref{eq:relB}), also because $j\circ g$ is totally geodesic.
	
	By the first formula in (\ref{eq:rel2}), with  $\bar f$, $j\circ g$ and $f$ playing the roles of $f$, $g$ and $F$ in that formula, respectively, and taking into account the first formula in Lemma \ref{le:comp} and the fact that $R^g=b^2 I$ for the totally geodesic embedding $g$, it follows that 
	$$\<T^f(j\circ g)_*\xi, (j\circ g)_*\eta\>=\<R^{j\circ g}\xi, \eta\>=\<R^g\xi, \eta\>=b^2\<\xi, \eta\>$$
	for all $\xi, \eta\in N_{\bar f}M$, hence $$\<(T-b^2I)(j\circ g)_*\xi, (j\circ g)_*\eta\>=0$$ for all $\xi, \eta\in N_{\bar f}M$.
	
Finally, to see that $\Phi=0$, recall from (\ref{eq:vartheta2}) that $\Phi=A^F_{\vartheta}$, where $F=h\circ f$. Write $h=k\circ \bar h$, where $\bar h\colon\, \Q_{k_1}^{n_1}\times \Q_{k_2}^{n_2}\to \Q_{k}^{n_1+n_2+1}$ and $k\colon\, \Q_{k}^{n_1+n_2+1}\to \R_{2\epsilon}^{n_1+n_2+2}$ are inclusions. Using that $\vartheta\in k_*N_{\bar h}(\Q_{k_1}^{n_1}\times \Q_{k_2}^{n_2})$ and that $\bar h\circ j\circ g$ is totally geodesic we obtain from (\ref{eq:relA}), with $\bar f$ and $h\circ j\circ g$ playing the roles of $f$ and $g$ in that formula, respectively, that
$$\< A^F_{\vartheta}X, Y\>=\<\alpha_{h\circ j\circ g}(\bar f_*X, \bar f_*Y), \vartheta\>=\<k_*\alpha_{\bar h\circ j\circ g}(\bar f_*X, \bar f_*Y), \vartheta\>=0.\qed\vspace{1ex}$$

	An isometric immersion $f\colon\, M^m\to N^n$ between Riemannian manifolds is said to be \emph{$1$-regular} if its first normal spaces have constant dimension on $M^m$.

 \begin{corollary}\po\label{cor:redcod1}  Let $f\colon\, M^m\to \Q_{k_1}^{n_1}\times \Q_{k_2}^{n_2}$, $k_1k_2\neq 0$,  be a $1$-regular  isometric immersion. Assume that $\Phi=0$ and that $N_1$ is parallel in the normal connection. Set $\ell=\rank N_1$.  Then $k_1k_2>0$, $n_i\geq m+\ell$ for $1\leq i\leq 2$, and there exist an isometric immersion $\bar f\colon\, M^m\to \Q_{k}^{m+\ell}$, with $k=k_1k_2/(k_1+k_2)$, and  a totally geodesic inclusion $j\colon\,\Q_{k_1}^{m+\ell}\times \Q_{k_2}^{m+\ell}\to \Q_{k_1}^{n_1}\times \Q_{k_2}^{n_2}$ such that $$f=j\circ g\circ \bar f,$$ where $g\colon\,  \Q_{k}^{m+\ell}\to \Q_{k_1}^{m+\ell}\times \Q_{k_2}^{m+\ell}$ is the totally geodesic embedding of Example \ref{ex:tgeod}.
 \end{corollary}
 \proof Since  $R=b^2 I$, with $b:=\sqrt{\frac{k_1}{k_1+k_2}}$, then  $S(TM)\subset N_1^\perp$ by  Lemma \ref{le:nablaR}. By the assumption that  $N_1$ (hence $N_1^\perp$) is a parallel subbundle of $N_fM$ with respect to the  normal connection, we have that 
 $(\nabla_X S)Y\in N_1^\perp$ for all $X, Y\in TM$. Then (\ref{eq:pid2}) implies that $$(T-b^2 I)N_1\subset N_1^\perp,$$ and the statement follows from Theorem \ref{thm:redcod2}.\qed

\section[Parallel submanifolds]{Parallel submanifolds}

We now use the results of the previous sections to prove Theorems \ref{thm:parallel} and \ref{thm:parallel2} in the introduction.

\begin{lemma}\po\label{le:SPhi} Let $f\colon\,M^m\to \Q_{k_1}^{n_1}\times \Q_{k_2}^{n_2}$ be a parallel  isometric immersion. Then the following holds:
\begin{itemize}
\item[$(i)$] If $k_1\neq 0$ and $k_2=0$, then at any point $x\in M^m$ either $S=0$ or there exist a unit vector $B\in T_xM$ and $\lambda\in (0,1)$ such that 
\be\label{eq:SR}
S X = \<X,B\>S B\,\,\,\,\mbox{and}\,\,\,(I-R)X=\lambda\<X,B\>B.\ee
\item[$(ii)$] If $k_1k_2 \neq 0$ then either $S$ or $\Phi$ vanishes everywhere on $M^m$. 
\end{itemize}
\end{lemma}
\proof Given $x\in M^m$, from  Codazzi equation (\ref{eq:codazzi})  we obtain
\be\label{eq:ncod}\<\Phi X, Z\>S Y = \<\Phi Y, Z\>S X, \,\,\,\mbox{for all}\,\,\, X, Y, Z \in T_x M.\ee
If $\Phi\neq 0$, then  either $S=0$ or the first equation in (\ref{eq:SR}) holds for 
 a unit vector $B\in T_xM$ spanning $(\ker S)^\perp$.
Replacing that equation into (\ref{eq:ncod}) yields
	$$\Phi X=\<X,B\>\Phi B=\mu\<X,B\>B$$
	for some $\mu\neq 0$ and for all $X\in T_xM,$ where in the last equality we have used that $\Phi$ is self adjoint. 
	
	Suppose  that  $k_2=0$. Then 
	$\Phi=k_1(I-R)$, hence $S=0$ wherever $\Phi=0$. Therefore  $(i)$ follows with $\lambda=\mu/k_1$. 
	
	Assume now that $k_1k_2\neq 0$. Then  we get a contradiction by assuming that  both $\Phi$ and $S$ are nonvanishing at $x$. In fact, in this case
	we would have that $\{B\}^\perp\subset  \ker \Phi\cap \ker S$, hence it suffices to prove that 
\be\label{eq:W} W:=\ker \Phi\cap \ker S=\{0\}\ee
if $k_1k_2\neq 0$. Otherwise we would have  $k_1I|_W-(k_1+k_2)R|_W=0$. Thus $k_1+k_2\neq 0$ and $R|_W=\frac{k_1}{k_1+k_2}I|_W$, hence the first equation in (\ref{eq:pi1}) would give
 $$		0 = S^t S|_W = R (I - R)|_W = \frac{k_1k_2}{\left(k_1 + k_2\right)^2}I|_W,$$
a contradiction.

	Now let $$A := \{x \in M^m \,:\, S|_{T_xM} \neq 0 \}\,\,\,\mbox{and}\,\,\,B:= \{x \in M^m \,:\, \Phi|_{T_xM} \neq 0\}.$$  Clearly, both $A$ and $B$ are open subsets of $M^m$. We have just proved that  $A^c \cup B^c = M^n$. On the other hand, it follows from (\ref{eq:W})  that $A^c \cap B^c = \emptyset$. Thus, either $A^c=M^n$ or $B^c=M^m$.\vspace{2ex}\qed
	
	\subsection{Proof of Theorem \ref{thm:parallel}.}
	
	 By Lemma \ref{le:SPhi}, either $S$ or $\Phi$ vanishes identically on $M^m$. 
Suppose first that $S$ is identically zero. If  either $R=0$ or $R=I$ everywhere, then case $(i)$ in the statement holds by  Proposition \ref{prop:R=0}. Otherwise,  we obtain from  Proposition \ref{prop:S=0} that $M^m$ is locally a Riemannian product $M_1 \times  M_2$ and $f=f_1\times f_2$, where $f_i\colon\, M_i \to \Q_{k_i}^{n_i}$, $1\leq i\leq 2$, is an isometric immersion. Since $f$ is parallel, it is easily seen that the same must hold for $f_i$, $1\leq i\leq 2$.

	Assume now that $\Phi = 0$. First observe that parallelism of the second fundamental form implies that $f$ is  $1$-regular and that $N_1$ is parallel in the normal connection. Then Corollary \ref{cor:redcod1}  applies. We obtain that  $k_1k_2>0$, $n_i\geq m+\ell$ for $1\leq i\leq 2$, where $\ell=\rank N_1$, and that there exist an isometric immersion $\bar f\colon\, M^m\to \Q_{k}^{m+\ell}$, with $k=k_1k_2/(k_1+k_2)$, and a totally geodesic inclusion $j\colon\,\Q_{k_1}^{m+\ell}\times \Q_{k_2}^{m+\ell}\to \Q_{k_1}^{n_1}\times \Q_{k_2}^{n_2}$ such that $f=j\circ g\circ \bar f,$ where $g\colon\,  \Q_{k}^{m+\ell}\to \Q_{k_1}^{m+\ell}\times \Q_{k_2}^{m+\ell}$ is the totally geodesic embedding of Example \ref{ex:tgeod}. Again, since $f$ is parallel, the same must hold for $\bar f$.\qed
 
 \subsection{The case $k_2=0$}
 
 The following reduction of codimension theorem for parallel submanifolds of symmetric spaces was obtained by Dombrowski \cite{dom}.
 
\begin{theorem}\po \label{thm:dom}  Let $N$ be a symmetric space. If $f\colon\,  M\to  N$ is a parallel isometric
immersion and if for some point $x\in M$ the second osculating space ${\cal O}_xf = f_*T_xM\oplus N_1(x)$ is contained in
some curvature invariant subspace $V$  of $T_{f(x)}N$, then $f(M)\subset \bar N$, where $\bar N$ denotes the totally
geodesic submanifold $\exp^N_{f(x)}(V)$.
\end{theorem}

We will make use of the following consequence of the preceding theorem.

\begin{corollary}\po \label{cor:dom} Let $f\colon\,  M\to  N$ be a parallel isometric
immersion into a  symmetric space. Assume that there exists an open subset ${\cal U}\subset M$ such that $f({\cal U})$ is  contained in
a totally geodesic submanifold $\bar N\subset N$. Then $f(M)\subset \bar N$.
\end{corollary}

\noindent \emph{Proof of Theorem \ref{thm:parallel2}:}
Assume first that $f\colon\,M^m\to \Q_{k_1}^{n_1}\times \R^{n_2}$, $k_1\neq 0$, is a  parallel isometric immersion such that  $S= 0$ everywhere.
If either $R=0$ or $R=I$, then $f$ is as in $(i)$  by Proposition \ref{prop:R=0}. Otherwise, it is given as in $(ii)$ by Proposition \ref{prop:S=0}.

Suppose now that $S\neq 0$ on an open subset ${\cal U}\subset M$. By Lemma \ref{le:SPhi}, the tensor $R$ is given on ${\cal U}$ by (\ref{eq:SR}) for some unit vector field $B$ and some smooth real function $\lambda$ with values in $(0,1)$.  Notice that the first equation in (\ref{eq:SR}) can also be written as
 $$S^t\eta=\<SB, \eta\>B$$
 for any $\eta\in N_fM$. In particular, $\ker S^t=\{SB\}^\perp$ splits orthogonally as $\ker S^t=U\oplus V$, with $U=\ker T$ and 
 $V=\ker (I-T)$. Arguing as in the proof of Lemma \ref{le:kerS}, both  $U$ and $V$ have constant rank.

 Equation (\ref{eq:pid2}) is equivalent to
 \be\label{eq:pid3} (\nabla_XS^t)\xi=A_{T\xi}X-RA_\xi X\ee
 for all $X\in TM$, $\xi\in N_fM$. For $\xi\in U=\ker T$ it yields
 \be\label{eq:sb}\<SB, \nabla_X^\perp \xi\>B=S^t\nabla^\perp_X \xi=-(\nabla_XS^t)\xi=RA_\xi X=A_\xi X-\lambda\<A_\xi X, B\>B\ee
 for all $X\in TM$. 
Therefore
\be\label{eq:rho}A_\xi X=\rho\<X, B\>B\ee
for some $\rho\in C^\infty({\cal U})$. 
In particular, if  $\xi\in U$ is orthogonal to $\zeta:=(\alpha(B, B))_U$, then $\xi\in N_1^\perp$.

   Now, given $X\in TM$ and $\xi\in N_fM$, we  have
   $$\pi_2(f_*X+\xi)=f_*(RX+S^t\xi)+SX+T\xi,$$
   hence $\pi_2(f_*X+\xi)=0$ if and only if 
   \be\label{eq:pi20}X-\lambda\<X, B\>B=RX=-S^t\xi=-\<SB, \xi\>B\ee
   and 
   \be\label{eq:pi20b}T\xi=-SX=-\<X,B\>SB.\ee
   Write $X=\mu B+Y$, with $\<Y,B\>=0$, and $\xi=\beta SB+\eta_U+\eta_V$, with $\eta_U\in U$ and $\eta_V\in V$. 
Then (\ref{eq:pi20}) becomes
$$\beta|SB|^2B=-Y+(\lambda-1)\mu B,$$
hence $Y=0$ and $\beta\lambda=-\mu$, using that $|SB|^2=\lambda(1-\lambda)$. On the other hand,  (\ref{eq:pi20b}) gives
$$\beta TSB+\eta_V=-\mu SB,$$
hence
$$\beta \lambda SB+\eta_V=-\mu SB,$$
where we have used  that $TSB=S(I-R)B=\lambda SB$. We obtain that $\eta_V=0$ and conclude that the kernel of $\pi_2$ is  the  subspace spanned by $U$ and the vector $\lambda f_*B-SB$.

If the vector field $\zeta=(\alpha(B, B))_U$ vanishes everywhere on ${\cal U}$, then $U\subset N_1^\perp$. Thus $U\cap N_1^\perp=U$ has constant rank on ${\cal U}$, and since $N_1$ is parallel in the normal connection,
condition (\ref{eq:condred}) is trivially satisfied. We claim that $f$ reduces codimension on the left by  $n_1-1$ on ${\cal U}$. To prove our claim, it is equivalent to show that if $f$ does not reduce codimension on the left then $n_1\leq 1$. By Corollary \ref{cor:redcod}, we must have $U=\{0\}$. Then the claim follows from that fact that the kernel of $\pi_2$ is the one-dimensional subspace spanned by  $\lambda f_*B-SB$, which implies that $n_2\geq n_1+n_2-1$. Therefore $f({\cal U})$ is contained in a totally geodesic submanifold
$\Q_{k_1}^{1}\times\R^{n_2}$ of $\Q_{k_1}^{n_1}\times\R^{n_2}$, and hence $f(M^m)\subset \Q_{k_1}^{1}\times\R^{n_2}$ by Corollary \ref{cor:dom}. We conclude that $f$ is as in $(iii)$. 

To finish the proof of Theorem \ref{thm:parallel2}, it remains to show that if there exists an open subset ${\cal U}\subset M$ where $S\neq 0$ and  the vector field $\zeta=(\alpha(B, B))_U$ is nowhere vanishing, then $f$ is given as in $(iv)$. This is the more delicate part of the proof. 

 First notice that in this case $U\cap N_1^\perp=U\cap \{\zeta\}^\perp$, hence we have again that  $U\cap N_1^\perp$ has constant rank on ${\cal U}$ and satisfies 
condition (\ref{eq:condred}), since $N_1$ is parallel. We now argue that  $f$ reduces codimension on the left by  $n_1-2$ on ${\cal U}$. In fact, assuming as before that $f$ does not reduce codimension on the left, we obtain from Corollary \ref{cor:redcod} that $U\cap \{\zeta\}^\perp=\{0\}$, that is, $U=\spa\{\zeta\}$. It follows that $n_1\leq 2$, because  the kernel of $\pi_2$ is now  spanned by  $\lambda f_*B-SB$ and $\zeta$, which implies that $n_2\geq n_1+n_2-2$. Therefore $f({\cal U})$ is contained in a totally geodesic submanifold
$\Q_{k_1}^{2}\times\R^{n_2}$ of $\Q_{k_1}^{n_1}\times\R^{n_2}$, and hence $f(M^m)\subset \Q_{k_1}^{2}\times\R^{n_2}$ by Corollary \ref{cor:dom}.

\begin{lemma}\po\label{le:main} Let $f\colon\,M^m\to \Q_{k}^{2}\times \R^{n}$, $k\neq 0$, be a  parallel isometric immersion with $S\neq 0$ everywhere that does not reduce codimension on the left. Then $M^m$ is locally a Riemannian product $M^m=\R\times N^{m-1}$ and $f=\gamma\times \tilde f$, where  $\gamma\colon\,\R\to \Q_{k}^{2}\times \R$ is a full extrinsic circle and $\tilde f\colon\,N^{m-1}\to \R^{n-1}$ is a parallel isometric immersion. 
\end{lemma}
\proof  By (\ref{eq:SR})  we have
 $$SX=\<X,B\>SB\,\,\,\,\mbox{and}\,\,\,\,(I-R)X=\lambda\<X,B\>B,$$
 where $B$ is a locally defined smooth unit vector field and $\lambda$ is a nonvanishing smooth function. \vspace{1ex}\\
 {\bf Claim $1$:} $\spa \{B\}$ and $\{B\}^\perp$ are totally geodesic distributions.\vspace{1ex}\\
 By  the assumption that $f$ does not reduce codimension on the left we have 
 $$N_fM=\spa\{\xi\}\oplus \spa\{SB\}\oplus V,$$
 where $V=\ker (I-T)$ and $\xi$ is a unit vector field spanning $U=\ker T$ for which (\ref{eq:rho}) holds with $\rho\neq 0$ on any open subset. From (\ref{eq:sb}) and (\ref{eq:rho}) we obtain 
 $$\<\nabla_X^\perp \xi, SB\>=\rho(1-\lambda)\<X, B\>.$$
 On the other hand, using that $T\xi=0$ and $S^t\xi=0$, it follows from (\ref{eq:pid4}) that 
 $$-T\nabla^\perp_X\xi=(\nabla^\perp_X T)\xi=-SA_\xi X-\alpha(X, S^t\xi)=-\rho\<X, B\>SB,$$
 hence
 $$\<\nabla^\perp_X\xi, \eta\>=0$$
 for any $\eta\in V$. It follows  that 
\be\label{eq:nperp}
\nabla^\perp_X\xi=|SB|^{-2}\<\nabla^\perp_X\xi, SB\>SB=\frac{\rho}{\lambda}\<X,B\>SB.
\ee
 By equation (\ref{eq:pid3}) for $\xi=SB$, i.e.,
  $$\nabla_XS^tSB-S^t\nabla_X^\perp SB=A_{TSB}X-RA_{SB}X,$$
  we obtain that
  $$\lambda\<\nabla_X Y, B\>=\<A_{SB}X, Y\>$$
  for all $X\in TM$ and $Y\in \{B\}^\perp$. In particular, 
  \be\label{eq:bb1}
  (A_{SB}B)_{\{B\}^\perp}=-\lambda\nabla_BB.
  \ee
  
  Since $f$ is parallel, we have
  $$\nabla_Y A_\xi X=A_\xi\nabla_Y X+A_{\nabla^\perp_Y\xi}X.$$
  Using (\ref{eq:rho}), the left-hand-side of the preceding equation becomes
  $$
  Y(\rho)\<X, B\>B+\rho\<\nabla_YX, B\>B+\rho\<X, \nabla_Y B\>B+\rho\<X, B\>\nabla_YB,$$
  whereas the right-hand-side is
  $$\rho\<\nabla_YX, B\>+\frac{\rho}{\lambda}\<Y, B\>A_{SB}X,$$
  where we have used (\ref{eq:nperp}). 
  Therefore,
  $$Y(\rho)\<X, B\>B+\rho\<X, \nabla_Y B\>B+\rho\<X, B\>\nabla_YB=\frac{\rho}{\lambda}\<Y, B\>A_{SB}X.$$
  For $\<X, B\>=0=\<Y, B\>$ we obtain
  $$\<\nabla_Y X, B\>=0.$$
  Hence $\{B\}^\perp$ is totally geodesic.
  For $\<X, B\>=0$ and $Y=B$ we get
  \be\label{eq:asb}A_{SB}X=\lambda\<\nabla_B B,X\>B,\ee
  which implies that 
  \be\label{eq:bb2}
  (A_{SB}B)_{\{B\}^\perp}=\lambda\nabla_BB.
  \ee
  Comparing (\ref{eq:bb1}) and (\ref{eq:bb2}) yields
  \be\label{eq:nbb}
  \nabla_B B=0,\ee
  and the proof of the Claim $1$ is completed.\vspace{1ex}
  
  Let $h\colon\, \Q_{k}^{2}\times \R^{n}\to \R_{\sigma(k)}^{n+3}$ be the inclusion and set $F=h\circ f$. \vspace{2ex}\\
    {\bf Claim $2$\,:} The second fundamental form of $F$ satisfies $\alpha_F(B, X)=0$ for any $X\in \{B\}^\perp$.\vspace{1ex}\\
  We have from (\ref{eq:sffj}) and (\ref{eq:relA}) that
 $$\alpha_F(X,Y)=h_*\alpha_f(X,Y)-k\<(I-R)X, Y\>\tilde \pi_1\circ F=h_*\alpha_f(X,Y)-k\lambda\<X, B\>\<Y, B\>\tilde \pi_1\circ F,$$
 hence it suffices to prove that $\alpha_f(B, X)=0$ for any $X\in \{B\}^\perp$, or equivalently, that $\{B\}^\perp$ is invariant under $A_\zeta$ for any $\zeta\in N_fM$. 
 
 We have from (\ref{eq:rho}) that $\{B\}^\perp$ is invariant under $A_\xi$,  whereas (\ref{eq:asb}) and (\ref{eq:nbb}) imply that
 the same holds for  $A_{SB}$.  Thus, it remains to show that this holds for any $\eta\in V$.

   By (\ref{eq:pid2}) and (\ref{eq:SR}) we have
   $$\nabla_X^\perp SB-S\nabla_XB=T\alpha(X, B)-\alpha(X, RB)=(T-I)\alpha(X, B)+\lambda\alpha(X, B).$$
   Taking the inner product with $\eta\in V=\ker (I-T)$ yields
  \be\label{eq:psi2}\<\nabla_X^\perp SB, \eta\>=\lambda\<A_\eta X, B\>.\ee
   On the other hand, since $f$ is parallel we have
   $$\nabla_Y A_{SB} X=A_{SB} \nabla_Y X+A_{\nabla_Y^\perp SB}X$$
   for all $X, Y\in TM$. Taking into account that $\{B\}^\perp \subset \ker A_{SB}$ by (\ref{eq:asb}) and (\ref{eq:nbb}), and that
   $\{B\}^\perp$ is totally geodesic, we obtain  
   \be\label{eq:psi1} A_{\nabla_Y^\perp SB}X=0\ee
   for all $X, Y\in \{B\}^\perp$. Define $\psi\colon\, \{B\}^\perp \to V$  by 
   $\psi(Y)=\nabla_Y^\perp SB.$ Then $\{B\}^\perp\subset \ker A_\eta$ for any $\eta\in \Ima \psi$ by (\ref{eq:psi1}), whereas (\ref{eq:psi2}) implies that $A_\eta(\{B\}^\perp)\subset \{B\}^\perp$ for any 
   $\eta$ in the orthogonal complement of $\Ima \psi$ in $V$.  Therefore $A_\eta\{B\}^\perp\subset \{B\}^\perp$ for any $\eta\in V$, and Claim $2$ is proved.

   It follows from Claim $1$ and the local version of de Rham theorem that $M^m$ splits locally as a Riemannian product $M^m=\R\times N^{m-1}$. Moreover, the splitting is global if
   $M^m$ is complete and simply connected by the global version of de Rham Theorem. 
   
   Define $$\R^{n-\ell}=\spa\{F_*(x)X\,:\,x\in M^m, \,X\in \{B\}^\perp(x)\}.$$
  Since $\{B\}^\perp\subset \ker (I-R)$ we have  that $\tilde \pi_1(F_*\{B\}^\perp)=\{0\}$, hence $\R^{n-\ell}\subset \R^{n}$. 
 By Claim~$2$ and   Moore's lemma \cite{mo}, 
there exist a  unit speed curve   $\gamma\colon\,\R\to \R^{\ell+3}=(\R^{n-\ell})^\perp$ and a full isometric immersion $\tilde{f}\colon\,N^{m-1}\to \R^{n-\ell}$ such that $F(x,y)=(\gamma(x), \tilde{f}(y))$.

Since  $F(M^m)\subset \Q_{k}^{2}\times \R^{n}$ and $f$ is parallel, then $\gamma(\R)\subset \Q_{k}^{2}\times \R^{\ell}$ and both $\tilde f$ and $\gamma$ are parallel.  It follows  from Corollary  \ref{cor:redcod} that $\gamma(\R)$ is contained in a totally geodesic $\Q_{k}^{2}\times \R\subset \Q_{k}^{2}\times \R^{\ell}$. Moreover, since $S\neq 0$ and $f$ does not reduce codimension on the left, $\gamma$ must be full in $\Q_{k}^{2}\times \R$. \vspace{1ex}\qed

We now complete the proof of Theorem \ref{thm:parallel2}. Let $f\colon\,M^m\to \Q_{k}^{2}\times \R^{n}$, $k\neq 0$, be a  parallel isometric immersion such that $S\neq 0$ and  the vector field $\zeta=(\alpha(B, B))_U$ is nowhere vanishing on an open subset ${\cal U}\subset M$. In view of Lemma \ref{le:main}, it is enough to show that $S$ can not vanish on any open subset.

Since the first factor in $\Q_{k}^{2}\times \R^{n}$ has dimension two, then $\rank(I-R)\in \{0,1,2\}$. If $R=I$ on an open subset of $M^m$, then $f$ is as in $(i)$ by Proposition \ref{prop:R=0} and Corollary~\ref{cor:dom}. Hence $R=I$ everywhere, contradicting the fact that $S\neq 0$ on ${\cal U}$. If $S= 0$ and $\rank(I-R)=2$ on a maximal open subset ${\cal V}$, then this also holds on $\bar{\cal V}$, since the eigenvalues of $(I-R)$ are either $0$ or $1$ on ${\cal V}$. But then $\rank(I-R)=2$ on an open neighborhood of $\bar{\cal V}$, contradicting the fact that $\rank(I-R)=1$ at points where $S\neq 0$ by Lemma \ref{le:SPhi}. 

 Therefore, if $S=0$ on a  maximal open subset ${\cal V}$,  there is a smooth unit vector field  $B$ such that $(I-R)X=\<X,B\>B$ for any $X\in T{\cal V}$. Let $x\in \bar{\cal V}$. Then there exist an open neighborhood $W$ of $x$ and  a smooth nowhere vanishing function $\lambda$ on $W$ such that $(I-R)X=\lambda\<X,B\>B$ for any $X\in TW$. The proof of Lemma \ref{le:main} shows that there exists a product neighborhood $I\times Z$ of $x$, an extrinsic circle $\gamma\colon\, I\to \Q^2_k\times \R$ and a parallel isometric immersion
$\bar f\colon\, Z\to \R^{n-1}$ such that $f|_{I\times Z}=\gamma\times \bar f$. Since $(I\times Z)\cap {\cal V}\neq \emptyset$, there exists an open interval $J\subset I$ such that  $\gamma(J)\subset \Q^2_k$, hence we must have $\gamma(I)\subset \Q^2_k$ by Corollary \ref{cor:dom}. But then $S=0$ on $W$, hence on ${\cal V}\cup W$, contradicting the maximality of ${\cal V}$ with respect to this property.\qed

\section{Umbilical submanifolds of $\Q_{k_1}^{n_1}\times \Q_{k_2}^{n_2}$}
 
 Let $f\colon\,M^m\to \Q_{k_1}^{n_1}\times \Q_{k_2}^{n_2}$ be an umbilical isometric immersion. Denote by $\eta$ its mean curvature vector. Then formulas (\ref{eq:pid1})--(\ref{eq:pid4}) become, respectively,
 \be\label{eq:nablaR}(\nabla_XR)Y=\<SY,\eta\>X+\<X,Y\>S^t\eta,\ee
	\be\label{eq:nablaS}(\nabla_XS)Y=\<X,Y\>T\eta-\<X,RY\>\eta\ee
	and
	\be\label{eq:nablaT}(\nabla_XT)\xi=-\<\xi,\eta\>SX-\<X,S^t\xi\>\eta.\ee
The Gauss, Codazzi and Ricci equations (\ref{eq:gauss}), (\ref{eq:codazzi}) and (\ref{eq:ricci}), respectively, 
 take the form
\be\label{eq:gaussumb} R(X,Y)=k_1(X\wedge Y-X\wedge RY-RX\wedge Y)+(k_1+k_2)RX\wedge RY+\|\eta\|^2X\wedge Y,\ee
\be\label{eq:codumb}\<Y,Z\>\nabla_X^\perp \eta-\<X,Z\>\nabla_Y^\perp \eta=\<\Phi X, Z\>SY-\<\Phi Y, Z\>SX\ee
and
\be\label{eq:ricciumb} R^\perp(X,Y)\xi=(k_1+k_2)(SX\wedge SY)\xi,\ee
whereas the equivalent form (\ref{eq:codazzi2}) of the Codazzi equation becomes
\be\label{eq:codumb2}\<\xi,\nabla_Y^\perp\eta\>X-\<\xi,\nabla_X^\perp\eta\>Y=\<S X, \xi\>\Phi Y-\<S Y, \xi\>\Phi X.\ee

\subsection{Case $S=0$}

Our approach to proving Theorem \ref{thm:umb} is based
on an analysis of the various possible  structures of the  tensor $S$. 
The next lemma takes care of the simplest case, in which $S$ vanishes identically.

\begin{lemma}\po\label{le:S=0} Let $f\colon\,M^m\to \Q_{k_1}^{n_1}\times \Q_{k_2}^{n_2}$ be an umbilical nontotally geodesic isometric immersion.  Assume that   $S=\{0\}$ everywhere. Then $f$ is as in part $(i)$ of Theorem \ref{thm:umb}.
\end{lemma}
\proof  Let $x\in M^m$ be a point where the mean curvature vector $\eta$ is nonzero  and let ${\cal U}\subset M^m$ be the maximal connected open neighborhood of $x$ where $\eta$ does not vanish.  It follows from (\ref{eq:nablaS}) that $R=\lambda I$ on ${\cal U}$, where $\lambda\in C^\infty({\cal U})$ is given by $\lambda=|\eta|^{-2}\<T\eta, \eta\>$. On the other hand, $R$ is an orthogonal projection at every $x\in {\cal U}$ by the first equation in (\ref{eq:pi1}), hence its eigenvalues are either $0$ or $1$. Thus, either $R=0$ or $R=I$ on ${\cal U}$. Assuming the first possibility, we conclude from Proposition \ref{prop:R=0} that $f|_{\cal U}$ is as in the statement.  In particular, $|\eta|$ is constant on ${\cal U}$. If ${\cal U}\neq M$, then $|\eta|$ would be nonzero  on some open subset containing ${\cal U}$, contradicting the maximality of ${\cal U}$ with respect to this property. \qed

\subsection{Case $\ker S=0$}

Our  next step is to  consider  the other extreme case, in which  $\ker S=\{0\}$ at some point.

\begin{lemma}\po\label{le:kerS=0} Let $f\colon\,M^m\to \Q_{k_1}^{n_1}\times \Q_{k_2}^{n_2}$ be an umbilical isometric immersion with $m\geq 3$ and  $k_1+k_2\neq 0$. Assume that   $\ker S=\{0\}$ at some point $x\in M^m$. Then  there exist umbilical isometric immersions $f_i\colon\, M^m\to \Q_{\tilde k_i}^{n_i}\subset \R^{n_i+1}$, $1\leq i\leq 2$, with $\tilde{k}_1=k_1\cos^2\theta$ and $\tilde k_2=k_2\sin^2\theta$ for some $\theta\in (0, \pi/2)$,  such that $f=(\cos \theta f_1,\sin\theta f_2)$. 
\end{lemma}
\proof Let ${\cal U}\subset M$ be the maximal  connected open  subset containing $x$ where   $\ker S=\{ 0\}$. Let $X_1, \ldots, X_m$ be an orthonormal diagonalizing frame for $R$, and let $\lambda_1, \ldots, \lambda_m$ be the corresponding eigenvalues. From Codazzi equation (\ref{eq:codumb}) we obtain 
$$\nabla^\perp_{X_i}\eta=((k_1+k_2)\lambda_j-k_1)SX_i$$
for all $1\leq i,j\leq m$ with $i\neq j$. Using that $m\geq 3$, the preceding equation implies that 
$$(k_1+k_2)\lambda_j=(k_1+k_2)\lambda_\ell$$
for all $1\leq j, \ell\leq m$. Since $k_1+k_2\neq 0$, it follows that there exists $\lambda\in C^\infty({\cal U})$ such that $\lambda_i=\lambda$ for all $1\leq i\leq m$.  From (\ref{eq:nablaR}) we obtain
$$X_i(\lambda)X_j+(\lambda I-R)\nabla_{X_i}X_j=\<SX_j, \eta\>X_i,$$
for all $1\leq i,j\leq m$ with $i\neq j$, hence $X_i(\lambda)=0$ for all $1\leq i\leq m$. Hence $R=\lambda I$ on ${\cal U}$, where $\lambda$ is a constant in $(0,1)$. But this implies that $R=\lambda I$ on the closure $\bar {\cal U}$ of ${\cal U}$, thus $\ker S=\{0\}$ also on $\bar {\cal U}$, and hence on an open neighborhood of ${\cal U}$, contradicting the maximality of ${\cal U}$ with respect to this property. We conclude that $R=\lambda I$ on $M$. The proof is completed by Proposition \ref{prop:rmid}.\qed

\subsection{Case $\dim \ker S=k\in (0, \dim M)$.}

   In this subsection we study umbilical isometric immersions $f\colon\,M^m\to \Q_{k_1}^{n_1}\times \Q_{k_2}^{n_2}$ for which the kernel of the tensor $S$  has constant dimension $k\in (0,m)$ on $M^m$.

\begin{lemma}\po\label{le:kerSneq0} Let $f\colon\,M^m\to \Q_{k_1}^{n_1}\times \Q_{k_2}^{n_2}$, $k_1+k_2\neq 0$,  be an umbilical isometric immersion. Assume that $\ker S$ has constant dimension $k\in (0,m)$ on $M^m$. Then $k=m-1$. Moreover, either $\ker S=\ker R$ on $M^m$ or $\ker S=\ker (I-R)$ on $M^m$. 
\end{lemma}
\proof Since $\ker S$ has constant dimension, the same holds for $\ker R$ and $\ker (I-R)$. Moreover, since $\ker S$ is nontrivial, the same must hold for either  $\ker R$ or $\ker (I-R)$. We will show that in the first case we must have  $k=m-1$ and $\ker S=\ker R$. By a similar argument one shows that ($k=m-1$ and) $\ker S=\ker (I-R)$ in the second case.

Applying (\ref{eq:codumb})  for $X\in \ker R$ and $Y=Z$ orthogonal to $X$ gives
\be\label{eq:neta}\nabla^\perp_X\eta=0\,\,\,\,\,\mbox{for all}\,\,\,X\in \ker R,\ee
 whereas for $Z=X\in \ker R$ and $Y\in (\ker R)^\perp\neq \{0\}$ it yields 
 \be\label{eq:neta1}\nabla^\perp_Y\eta=-k_1SY\,\,\,\,\,\mbox{for all}\,\,\,Y\in (\ker R)^\perp.\ee 

Choose $Y\in (\ker S)^\perp$. Then, applying (\ref{eq:codumb2}) for $\xi =SY$ and using that $k_1+k_2\neq 0$ we obtain for any
$X\in (\ker R)^\perp$ that
$$RX=\frac{\<SX, SY\>}{|SY|^2}RY.$$
Thus $R$ has rank one, and the remaining of the statement follows from $\ker R\subset \ker S$.\vspace{1ex}\qed

In the next result  we assume that $\ker S$ has constant dimension $m-1$ and that  $\ker S=\ker R$ everywhere, the case 
$\ker S=\ker (I-R)$ being completely similar. 

\begin{lemma}\po\label{le:rcodumb} Let $f\colon\,M^m\to \Q_{k_1}^{n_1}\times \Q_{k_2}^{n_2}$  be an umbilical isometric immersion such that $\ker S$ has constant dimension $m-1$ and   $\ker S=\ker R$ everywhere. Then $f$ reduces codimension on the left by either $n_1-m$ or $n_1-m-1$ and on the right by $n_2-1$.
\end{lemma}
\proof  To prove the statement, it is equivalent to show that if $f$  reduces codimension neither on the left nor on the right 
then $n_2=1$ and $n_1$ is either $m$ or $m+1$.

By the assumption, there exist a locally defined smooth unit vector field $B$   and a smooth nowhere vanishing function $\lambda$ such that
$$SX=\<X, B\>SB\,\,\,\,\,\mbox{and}\,\,\,\,\,RX=\lambda\<X,B\>B$$
for all $X\in TM$. Hence $S^t\xi=\<SB, \xi\>B$ for all $\xi \in N_fM$. We obtain from (\ref{eq:nablaR}) that 
\be\label{eq:nr1} X(\lambda)\<Y,B\>B+\lambda\<Y, \nabla_XB\>B+\lambda\<Y, B\>\nabla_XB=\<SB, \eta\>(\<Y, B\>X+\<X, Y\>B).\ee
For $X=Y=B$, the preceding equation gives $B(\lambda)=2\<SB, \eta\>$ and $\nabla_BB=0$. Hence $(\ker R)^\perp=\spa\{B\}$ is a totally 
geodesic distribution. Applying (\ref{eq:nr1}) for $Y\in \ker R=\{B\}^\perp$ yields
$$\<\nabla_X Y, B\>=-\frac{\<SB, \eta\>}{\lambda}\<X, Y\>.$$
Thus $\ker R$ is an umbilical distribution with mean curvature vector field $\varphi=\mu B$, where $\mu=-\<SB, \eta\>/\lambda$, i.e., 
$$(\nabla_XY)_{(\ker R)^\perp}=\<X, Y\>\varphi\,\,\,\,\,\mbox{for all}\,\,\,\,X, Y\in \ker R.$$
Applying (\ref{eq:nablaS}) for a unit vector field $Y=X\in \{B\}^\perp$ gives
\be\label{eq:teta} T\eta=-\<\nabla_X X, B\>SB=-\mu SB.\ee
Then, for $Y=B$ it yields
\be\label{eq:nsb} \nabla_X^\perp SB=-\<X, B\>(\mu SB+\lambda \eta)\ee
for any $X\in TM$. In particular,
\be\label{eq:nsb2} \nabla_X^\perp SB=0\,\,\,\,\mbox{for any}\,\,\,\,X\in \{B\}^\perp.\ee
On the other hand, for $X\in \{B\}^\perp$ and $Y=B$ we obtain from (\ref{eq:nr1}) that 
\be\label{eq:xlambda} X(\lambda)=0\,\,\,\,\mbox{for any}\,\,\,\,X\in \{B\}^\perp.\ee
It follows from (\ref{eq:neta}), (\ref{eq:nsb2}) and (\ref{eq:xlambda}) that
$$X(\mu)=0\,\,\,\,\mbox{for any}\,\,\,\,X\in \{B\}^\perp,$$
hence $\ker R$ is a spherical distribution, i.e., $\nabla_X \varphi \in \ker R$ for all $X\in \ker R$.

By the local version of Hiepko's theorem \cite{hi}, there exists locally an isometry 
$$\psi\colon\, I \times_\rho N^{m-1}\to M^m$$ from a 
warping product manifold, where $I\subset\R$ is an open interval and $\rho\in C^\infty(I)$ is the  warping function, that maps the leaves of the product foliation induced by the factors $I$ and $N^{m-1}$ into the leaves of $\spa\{B\}$ and $\{B\}^\perp$, respectively. 

Our next step is to prove that the vector fields $\eta$ and $SB$ are either linearly independent everywhere or 
linearly dependent everywhere. For that we consider the function $h\colon\, M^m\to \R$ given by
$$h(x)=|\eta(x)|^2|SB(x)|^2-\<\eta(x), SB(x)\>^2,$$
which vanishes precisely at the points where $\eta$ and $SB$ are  linearly dependent. Using (\ref{eq:neta}), (\ref{eq:neta1}) and (\ref{eq:nsb}), a straightforward computation gives
\be\label{eq:xh} X(h)=-2\<X, B\>\mu h\,\,\,\,\,\mbox{for all}\,\,\,\,X\in TM.\ee
 Therefore $\tilde h=h\circ \psi$ depends only on $I$ and we can write  (\ref{eq:xh}) as  $\tilde h'(t)=-2\mu(t)\tilde h(t)$, where we also write $\mu$ for $\mu\circ \psi$. Hence  $$\tilde h(t)=\tilde h(t_0)\exp \left(-2\int_{t_0}^t\mu(s)\,ds\right)$$ for any fixed $t_0\in I$. We obtain that each point of $M^m$ has an open neighborhood $V$ where $h$ is either identically zero or is nowhere vanishing. Therefore, the set $h^{-1}(0)$ is both open and closed, hence $h$ is either identically zero or  nowhere vanishing in $M^m$.

Let $U=\ker T$ and $V=\ker (I-T)$. We  prove that $\eta \in U\oplus \spa \{SB\}$. 
We have
$$\eta=\eta_U+\eta_V+\frac{\<\eta, SB\>}{|SB|^2}SB=\eta_U+\eta_V+\frac{\<\eta, SB\>}{\lambda(1-\lambda)}SB,$$
where $\eta_U$ and $\eta_V$ are the components of $\eta$ in $U$ and $V$, respectively. Then, using that $TSB=S(I-R)B=(1-\lambda)SB$ we obtain
$$T\eta=\eta_V+\frac{\<\eta, SB\>}{\lambda}SB=\eta_V-\mu SB.$$
Comparing with (\ref{eq:teta}) gives $\eta_V=0$, as we claimed.

Given $\zeta\in U\cap \{\eta\}^\perp=U\cap N_1^\perp$, it follows from (\ref{eq:neta}) and (\ref{eq:neta1}) that
$$\<\nabla_X^\perp \zeta, \eta\>=-\<\zeta, \nabla_X^\perp \eta\>=0,$$
hence $\nabla_X^\perp \zeta\in \{\eta\}^\perp$. Similarly, we have  $\nabla_X^\perp \zeta\in \{\eta\}^\perp$ for any $\zeta\in V=V\cap \{\eta\}^\perp=V\cap N_1^\perp$.

Since we are assuming that $f$ reduces codimension neither on the left nor on the right, it follows from Corollary \ref{cor:redcod} that $V=\{0\}=U\cap \{\eta\}^\perp$. 
Hence, the codimension  $n_1+n_2-m$ is either $1$ or $2$,  according as $\eta_U$ is zero or not.

Now,
\begin{eqnarray*}\pi_1\eta&=&-f_*S^t\eta+(I-T)\eta=-\<\eta, SB\>f_*B+\eta_U+\frac{\<\eta, SB\>}{\lambda(1-\lambda)}SB-\frac{\<\eta, SB\>}{\lambda}SB\\&=&-\frac{\<\eta, SB\>}{1-\lambda}((1-\lambda)f_*B-SB)+\eta_U.
\end{eqnarray*}
On the other hand, 
$$\pi_1f_*X=f_*(I-R)X-SX=f_*X$$
if $X\in \ker R$, and 
$$\pi_1f_*B=f_*(I-R)B-SB=(1-\lambda)f_*B-SB.$$
It follows that
$$\pi_1(f_*TM\oplus \spa \{\eta\})=f_*\ker R\oplus \spa\{(1-\lambda)f_*B-SB\}\oplus \spa\{\eta_U\},$$
hence $n_1$ is either $m$ or $m+1$, according as $\eta_U$ is zero or not.\qed

\subsection{Proof of Theorem \ref{thm:umb}:}  We use the following consequence of part $(ii)$ of Proposition $1$ of \cite{al}:

\begin{proposition}\po \label{prop:al} Let $N^n$, $n\geq 3$, be a Riemannian manifold, let $x\in N^n$, let $L^\ell$ be an $\ell$-dimensional subspace of $T_xN^n$, $2\leq \ell<n$, and let $H^*$ be a vector in $T_xN^n$ orthogonal to $L^\ell$. Assume that there exist a totally umbilical isometric immersion $f\colon\,\bar{N}^\ell\to N^n$ of a connected Riemannian manifold $\bar{N}^\ell$ and $\bar x\in \bar{N}^\ell$ such that $f(\bar x)=x$, $f_*T_{\bar x}\bar N^\ell=L^\ell$ and $H^f(\bar x)=H^*$. Suppose also that there exist a complete Riemannian manifold $\tilde{N}^{\ell+1}$, a totally geodesic isometric immersion $g\colon\,\tilde{N}^{\ell+1}\to N^n$ and $\tilde x\in \tilde{N}^{\ell+1}$ such that $g(\tilde x)=x$ and $g_*T_{\tilde x}\tilde{N}^{\ell+1}=L^\ell\oplus \spa\{H^*\}$. Then $f(\bar{N}^\ell)\subset g(\tilde{N}^{\ell+1})$.
\end{proposition}

 Let $f\colon\, M^m\to  \Q_{k_1}^{n_1}\times \Q_{k_2}^{n_2}$ be as in  Theorem \ref{thm:umb}. Since $f$ is not totally geodesic, it follows from  Proposition \ref{prop:al} that
  there exists no open subset of $M^n$ where the mean curvature vector $\eta$ of $f$ vanishes. 

 If $S$ vanishes everywhere on  $M^n$, then $f$ is as in either $(i)$  by Lemma \ref{le:S=0}. If $\ker S=\{0\}$ at some point $x\in M^m$, then $f$ is as in $(ii)$ by Lemma \ref{le:kerS=0}. Then, we can assume that there exists an open subset ${\cal U} \subset M^m$ where $\ker S$ has constant dimension $k\in (0, m)$ and  $\eta$ is nowhere vanishing. By Lemma \ref{le:kerSneq0}, we must have   $k=m-1$. Moreover, either $\ker S=\ker R$ on ${\cal U}$ or $\ker S=\ker (I-R)$ on ${\cal U}$. Assume, say, that the first  possibility holds. Then, Lemma \ref{le:rcodumb} implies that,   on ${\cal U}$,  $f$ reduces codimension on the left by either $n_1-m$ or $n_1-m-1$ and on the right by $n_2-1$. It follows from  Proposition \ref{prop:al} that the same must hold everywhere on $M^m$. Since $S\neq 0$ and $\eta\neq 0$ on ${\cal U}$, the codimension of $f$ can not be further reduced either on the left  or on the right. Hence $f$ is as in $(iii)$. Similarly,  if $\ker S=\ker (I-R)$ on ${\cal U}$,  we conclude that $f$ is as in $(iii)$ after interchanging  the factors.\qed

{\renewcommand{\baselinestretch}{1}

\hspace*{-20ex}\begin{tabbing} \indent\= Universidade Estadual de Londrina
\indent\indent\=  Universidade Federal de S\~ao Carlos\\
\> Rodovia Celso Garcia Cid km 380
\> Rodovia Washington Luiz km 235\\
\> 86051-980  -- Londrina -- Brazil  \>
13565-905 -- S\~{a}o Carlos -- Brazil \\
\> E-mail: brunomrs@uel.br \> E-mail: tojeiro@dm.ufscar.br
\end{tabbing}}

\end{document}